\documentclass[aihp]{imsart}

\RequirePackage{amsthm,amsmath,amsfonts,amssymb}
\RequirePackage[numbers]{natbib}
\RequirePackage[colorlinks,citecolor=blue,urlcolor=blue]{hyperref}
\RequirePackage{graphicx}

\usepackage{calc}
\newsavebox\CBox
\newcommand\hcancel[2][0.5pt]{%
  \ifmmode\sbox\CBox{$#2$}\else\sbox\CBox{#2}\fi%
  \makebox[0pt][l]{\usebox\CBox}%
  \rule[0.5\ht\CBox-#1/2]{\wd\CBox}{#1}}
\usepackage{blindtext}
\usepackage{soul}
\usepackage{amsmath,amsthm,amssymb,amscd,amsfonts,amsbsy}
\usepackage{color}
\usepackage{hyperref,url}
\usepackage{float}
\usepackage{hyperref}
\usepackage{enumerate}
\usepackage{enumitem}   
\usepackage{bbm}
\usepackage{cancel}
\usepackage{mathrsfs}
\usepackage{colonequals}

\usepackage{amsmath,amsfonts,amsthm,bm}
\usepackage{mathrsfs}  
\usepackage{xcolor}
\usepackage{amsfonts}
\usepackage{amssymb}
\usepackage{amsmath}
\usepackage{mathtools}
\usepackage{mathrsfs}  
\usepackage[normalem]{ulem}

\startlocaldefs
\numberwithin{equation}{section}

\theoremstyle{remark}
\newtheorem{definition}{Definition}[section]
\newtheorem{remark}[definition]{Remark}
\newtheorem{notation}[definition]{Notation}
\newtheorem{example}[definition]{Example}

\theoremstyle{plain}

\newtheorem{theorem}[definition]{Theorem}
\newtheorem{hypothesis}{Hypothesis}
\newtheorem{lemma}[definition]{Lemma}
\newtheorem{proposition}[definition]{Proposition}
\newtheorem{corollary}[definition]{Corollary}

\newtheorem {mtheorem} {\bf Theorem}

\newtheorem{manualtheoreminner}{\bf Claim}

\newtheorem{manualtheoreminnerb}{\bf Step}
\newenvironment{step}[1]{%
  \renewcommand\themanualtheoreminnerb{#1}%
  \manualtheoreminnerb
}{\endmanualtheoreminner} 

\def\R{\mathbb R}
\def\d{\mathrm d}
\def \d {\mathrm{d}}

\def \Pb {\mathbb{P}}

\newcommand{\Bigwedge}{\textstyle{\bigwedge}}

\usepackage{mathtools}

\newcommand{\revision}{}

\endlocaldefs

\begin{document}

\begin{frontmatter}
\title{The conditioned Lyapunov spectrum for random dynamical systems}
\runtitle{The conditioned Lyapunov spectrum for random dynamical systems}

\begin{aug}
\author[A]{\fnms{Matheus M.}~\snm{Castro}\ead[label=e1]{m.de-castro20@imperial.ac.uk}\orcid{0000-0002-2513-2830}},
\author[B]{\fnms{Dennis}~\snm{Chemnitz}\ead[label=e2]{dennis.chemnitz@fu-berlin.de}\orcid{0000-0002-3303-3533}},
\author[A]{\fnms{Hugo}~\snm{Chu}\ead[label=e3]{hugo.chu17@imperial.ac.uk}\orcid{0000-0002-1620-3334}},
\author[B, KdV]{\fnms{Maximilian}~\snm{Engel}\ead[label=e4]{maximilian.engel@fu-berlin.de}\orcid{0000-0002-1406-8052}}
\author[A,IRCN,CAMB]{\fnms{Jeroen S.W.}~\snm{Lamb}\ead[label=e5]{jsw.lamb@imperial.ac.uk}\orcid{0000-0001-7647-4200}}
\and
\author[A]{\fnms{Martin}~\snm{Rasmussen}\ead[label=e6]{m.rasmussen@imperial.ac.uk}\orcid{0000-0002-7366-4719}}
\address[A]{Department of Mathematics,
Imperial College London, London SW7 2AZ, UK\printead[presep={,\ }]{e1,e3,e5,e6}}

\address[B]{Department of Mathematics, Freie Universität Berlin, Arnimallee 3,
14195 Berlin, Germany\printead[presep={,\ }]{e2,e4}}

\address[IRCN]{International Research Center for Neurointelligence, The University of Tokyo, Tokyo,113-0033, Japan}
\address[CAMB]{Centre for Applied Mathematics and Bioinformatics, Department of Mathematics and
Natural Sciences, Gulf University for Science and Technology, Halwally 32093, Kuwait}

\address[KdV]{Korteweg-de Vries Institute for Mathematics, University of Amsterdam,
1098 XG Amsterdam, The Netherlands}
\end{aug}

\begin{abstract}
We establish the existence of a full spectrum of Lyapunov exponents for memoryless random dynamical systems with absorption. To this end, we crucially embed the process conditioned to never being absorbed, the $Q$-process, into the framework of random dynamical systems, allowing us to study multiplicative ergodic properties. We show that the finite-time Lyapunov exponents converge in conditioned probability and apply our results to iterated function systems and stochastic differential equations.
\end{abstract}

\begin{abstract}[language=french]
Nous établissons l'existence d'un spectre complet d'exposants de Lyapunov pour les systèmes dynamiques aléatoires sans mémoire avec absorption. Pour celà, nous adaptons le processus conditionné à ne jamais être absorbé à la structure des systèmes dynamiques aléatoires, nous permettant ainsi d'étudier ses propriétés multiplicatives ergodiques. Nous montrons que la convergence vers les exposants de Lyapunov se produit en probabilité conditionnelle et nous appliquons nos résultats aux systèmes de fonctions itérées et aux équations différentielles stochastiques.
\end{abstract}



\end{frontmatter}


\section{Introduction}

A central part of modern mathematical theory and modelling is the description of evolving systems subject to uncertainty. A classical object of study are Markov processes on some state space $E$ which are given by a tuple $\left(\Omega, (\mathcal{G}_t), (X_t), (\mathcal P^t), (\mathbb{P}_x) \right)$. Here, the law of the $\mathcal{G}_t$-adapted stochastic process $X_t$ under the probability measure $\mathbb P_x$ describes the evolution of the modelled system, initialised at $x \in E$, giving rise to a semigroup structure $\mathcal P^t$. In that sense, this formalism only describes the statistics of the one-point motion of trajectories and joint probability distributions for different initial conditions are not defined. Hence, classical questions from dynamical systems, in particular concerning the sensitivity on initial conditions associated to chaos, cannot be addressed. The correct framework for studying such questions is given by the theory of random dynamical systems (RDS) \cite{Arnold1998RandomSystems}  which model the stochastic system as a (deterministic) skew product $(\theta, \varphi)$ where $\varphi$ evolves as a cocycle over the underlying noise dynamics given by $\theta$. In fact, since every RDS with independent increments induces a Markov process in a canonical way, it contains, in principle, more information.
Specifically, the framework of RDS allows for the definition and analysis of Lyapunov exponents which describe the asymptotics of the sensitivity to initial conditions.
However, for systems with a unique ergodic component, such as for those driven by unbounded noise, the classical theory of Lyapunov exponents only captures global dynamical properties, for instance, the  contraction of bounded sets to a single random fixed point \cite{Crauel1998AdditiveBifurcation,  Flandoli2017SynchronizationNoise,Flandoli2017SynchronizationSystems}.
This is one of the reasons why a stochastic extension of the local bifurcation theory for deterministic dynamical systems, describing changes of stability in dynamical behaviour, has been only developed along single examples and phenomena \cite{Arnold1999TheBifurcation, Arnold1996TowardStudy, Baxendale1994ABifurcation, Doan2018HopfNoise, Engel2019BifurcationCycle}.

{\revision A first step towards a description of local stability properties in random dynamical systems was undertaken by Engel, Lamb and Rasmussen \cite{Engel2019ConditionedSystems} in the context of stochastic differential equations (SDEs) with additive noise, extending the notion of a (dominant) Lyapunov exponent to dynamics conditioned to remain within a bounded subdomain of the state space. 

This \emph{conditioned Lyapuonov exponent} has already shown to be very useful in practice. Breden and Engel \cite{Breden2021Computer-assistedSystems} used rigorous computation to prove the existence of a noise-induced transition from negative to positive conditioned Lyapunov exponent, establishing in an adapted model of shear-induced chaos \cite{LinYoung2008} the existence of a transition from local noise-induced synchronisation to chaos (cf.~also \cite{ChemnitzEngel2023} for a result in the global  setting)}.
Bassols-Cornudella and Lamb \cite{bassolscornudella2023noiseinduced} have exploited conditioned Lyapunov exponents to reveal the mechanism behind a noise-induced transition to chaos in a random logistic map, modelling the interaction between effectively expanding and contracting compartments.

In this paper we extend the results of \cite{Engel2019ConditionedSystems}, establishing the existence of a full Lyapunov spectrum with corresponding Oseledets spaces. We overcome  limitations in \cite{Engel2019ConditionedSystems} due to its reliance on a conditioned version of the Furstenberg--Khasminskii formula for additive noise, with strong assumptions on the projective bundle process.
Instead, we provide a more appropriate, general framework for addressing dynamical questions in a conditioned setting, by  translating the notion of stationarity for asymptotic survival processes to a suitable invariant measure for the conditioned RDS. Specifically, this allows us to apply the Multiplicative Ergodic Theorem to a large class of conditioned stochastic processes, yielding in particular a full Lyapunov spectrum under relatively mild assumptions.

In more detail, consider a RDS $(\theta, \varphi)$ {\revision in one-sided time $\mathbb T = \mathbb N_0 := \mathbb N \cup \{0\}$ or $\mathbb T = \mathbb R_+ := [0,\infty)$} with filtered {\revision memoryless } probability space $(\Omega, (\mathcal{F}_t)_{t\in \mathbb T}, \mathcal{F}, \mathbb P)$   on a state space $E_M$ which is decomposable as $E_M=M\sqcup \{\partial\}$; here $M$ is a Riemannian manifold and ${\partial}$ is a cemetery (or absorbing) state for $\varphi$, i.e.~$\varphi_s\in \{\partial\}$ implies that $\varphi_t\in \{\partial \}$ for all $t\geq s$.\footnote{ {\revision See Appendix \ref{RDS} for a more detailed description of our setting.}} Accordingly, we introduce the stopping time
$$\tau(\omega, x):= \inf_{t\in \mathbb T}\left\{\varphi_t (\omega, x) = \partial\right\},\qquad (\omega, x) \in \Omega \times M.$$
The two classes of RDSs we consider are those given by solutions of SDEs (in continuous time) and iterations of random maps (in discrete time).
As indicated above, $(\theta, \varphi)$ induces a  Markov process 
$$\varphi:=\left(\Omega \times E_M, (\mathcal G_t:= \mathcal{F}_t \otimes \mathcal{B}(E_M))_{t\in \mathbb T}, (\varphi_t)_{t\in \mathbb T}, (\mathcal P^t)_{t\in \mathbb T}, (\mathbb{P}_x) _{x\in E_M}\right),$$
where $\mathbb P_x  := \mathbb P \otimes \delta_x$ and the usage of $\varphi$ as cocycle or Markov process becomes clear from the context.
Conditioning a stochastic system to never reach the cemetery state is a well-studied problem for Markov processes~\cite{Collet2013Quasi-StationaryDistributions, Pinsky1985AMeasures, Pinsky1985OnProcesses}, going back to the pioneering work of Yaglom \cite{Yaglom1947CertainProcesses}, with recent advances~\cite{Castro2021ExistenceApproach, Champagnat2016ExponentialQ-process, Champagnat2017UniformQ-process, Colonius2021Quasi-ergodicChains} on the statistical properties of the conditioned process, in particular on its ergodic properties. {\revision{This literature provides readily verifiable assumptions for our main hypothesis below which amounts to the exponential convergence of the statistics of the conditioned process to a \emph{quasi-stationary} distribution.}}

\begingroup
\renewcommand\thehypothesis{(H)} 
\revision{
\begin{hypothesis}\label{(H)} Let $(\Omega, (\mathcal{F}_t)_{t\in \mathbb{T}}, \mathcal{F}, \Pb, (\theta_t)_{t\in \mathbb{T}})$ be a memoryless noise space of the  form \eqref{eqn:cts} or \eqref{eqn:discrete}, and $(\theta,\varphi)$ a
random dynamical system on $E_M$ absorbed at $\partial$.
\begin{enumerate}
    \item[${\mathrm{(H1)}}$] The Markov process $(\Omega \times E_M, (\mathcal G_t)_{t\in \mathbb{T}}, (\varphi_t)_{t\in\mathbb{T}}, (\mathcal P^t)_{t\in \mathbb T}, (\mathbb{P}_x )_{x\in E_M})$ admits a unique quasi-stationary distribution $\mu$.
    \item[${\mathrm{(H2)}}$] There exists a constant $\alpha>0$ such that for every $x\in M$, there exist $C(x)$ such that
    \begin{align}
        \|\mathbb P_x(\varphi_t \in \cdot \mid  \tau > t) - \mu\|_{TV} \leq C(x) e^{-\alpha t}. \label{lim1}
    \end{align}
    \item[${\mathrm{(H3)}}$]  There exists a positive bounded function $\eta$ on $M$  and a constant $\beta >0,$ such that
    \begin{equation}\lim_{t\to\infty}\sup_{x\in M}\left| e^{\beta  t}\mathbb P_x(\tau > t) -\eta(x)\right|  = 0. \label{eqn:etalim}
    \end{equation}
\end{enumerate}
\end{hypothesis}}
\endgroup

A key ingredient for the following is the notion of the $Q$-process\cite{Champagnat2016ExponentialQ-process, Champagnat2017UniformQ-process}, which describes the process $\varphi$ conditioned on asymptotic survival. It is given by the $Q$-measures
$$\mathbb{Q}_{x}(A) := \lim_{t\to\infty}\mathbb{P}_x(A | \tau>t)\qquad \text{for all } A \in \mathcal{F}_s \otimes \mathcal{B}(E_M) , \text{ for any fixed } s > 0. $$
In the setting of Hypothesis~\ref{(H)}, it can be shown that these limits exist, these measures define a Markov process and the measure given by $\nu (\d x) = \eta(x) \mu(\d x)$ is a stationary distribution of this process. Sometimes this measure is also called quasi-ergodic \cite{Breyer1999AProcesses} (See Definiteion \ref{Def:QEM}) since it turns out that the Birkhoff averages satisfy
$$\lim_{t\to\infty}\mathbb{E}_x\left[\frac{1}{t}\int_0^t f(\varphi_t)\d t \bigg|~ \tau>t\right] = \int_M f\d \nu \qquad \text{for all } f \in L^1(\nu).$$
In \cite{Engel2019ConditionedSystems}, this property was exploited to obtain the notion of a dominant conditioned Lyapunov exponent via a modified Furstenberg--Khasminskii formula. Specifically, it was shown that for additive noise SDEs with linearisation $\mathrm D\varphi_t$, the following limit exists 
\begin{equation*}
    \Lambda_1 = \lim_{t\to\infty}\mathbb{E}_x\left[ \frac{1}{t}\log \|\mathrm{D}\varphi_t \| \,\bigg|~\, \tau>t\right].
\end{equation*}
Consequently, it was conjectured in \cite[Conjecture 3.5]{Engel2019ConditionedSystems} that additional exponents $\left\{\Lambda_i\right\}_{i=1}^d$ can be found as limits
\begin{equation}\label{eqn:defLyap}
\Lambda_i = \lim_{t\to\infty}\mathbb{E}_x\left[\frac{1}{t}\log\delta_i(\mathrm D\varphi_t)\,\bigg|\, \tau>t\right] \qquad i \in \{ 1, \ldots, d\},
\end{equation}
where $\delta_i(\mathrm D\varphi_t)$ denotes the $i\textsuperscript{th}$ singular value of $\mathrm D\varphi_t$.
To show this conjectured existence of a spectrum of conditioned Lyapunov exponents, we now find an appropriate invariant, ergodic measure for the random dynamical system corresponding to the quasi-ergodic distribution.
\begin{mtheorem}[Ergodic measure for conditioned RDS]\label{thm:IntroErgMeas}
    Let $\Theta := (\theta, \varphi)$ be a random dynamical system on $M$ with absorption at $\{\partial\}$ satisfying Hypothesis \ref{(H)} with quasi-ergodic distribution $\nu$.
    Then  $\Theta$ has an invariant, ergodic (even strongly mixing) probability measure given by
        $$\mathbb Q_{\nu}(\cdot) := \int_M \mathbb Q_x(\cdot) \nu(\d x).$$
\end{mtheorem}
This new crucial insight allows for the application of the multiplicative ergodic theorem to obtain the following theorem as a corollary:

\begin{mtheorem}[Lyapunov spectrum for the $Q$-process]\label{thm:IntroLyapSpect}
   Assume that the linear cocycle $\Phi := (\mathrm{D}\varphi_t)_{t\geq 0}$, as the linearisation over the $\mathcal C^1$ random dynamical system $\Theta := (\theta, \varphi)$ as in Theorem~\ref{thm:IntroErgMeas}, is invertible and fulfills the integrability condition
        \begin{equation}\label{eqn:IntroICMET}
            \mathbb{E}^{\mathbb{Q}}_{\nu}\left[\sup_{0\leq t\leq 1}\log^{+}\left\|\Phi^{\pm1}_t\right\|\right]<\infty.
        \end{equation}
        Then there exists a full spectrum of constant Lyapunov exponents $\Lambda_1\geq \cdots \geq \Lambda_d > -\infty$ such that for all $i\leq d$
        $$\lim_{t\to\infty}\mathbb{E}^{\mathbb{Q}}_{\nu}\left[\left|\Lambda_i- \frac{1}{t}\log\delta_i(\Phi_t)\right|\right] = 0.$$
\end{mtheorem}
Here, the expression $\mathbb{E}^{\mathbb{Q}}_{\nu}$ denotes expectation with respect to the measure $\mathbb{Q}_{\nu}$.
In more detail, we also obtain Oseledets flags, i.e.~dynamically invariant subspaces that constitute a filtration of the tangent space, which are associated with the distinct Lyapunov exponents (cf.~Theorem~\ref{thm:fullMET}). We remark that condition~\eqref{eqn:IntroICMET} can be verified via the following similar assumption in terms of the QSD $\mu$ 
$$ \mathbb E_\mu \left[\sup_{0\leq t\leq 1}\log^+ \|\Phi^{\pm1}_t \| \mathbbm 1_{\{\tau >1\} } \right] < \infty,$$
which is, in general, easier to verify.
Finally, we use results for the $Q$-process to show convergence to Lyapunov exponents in conditional probability, and under stronger assumptions, that are satisfied for SDEs with additive noise, even convergence in conditional expectation. In particular, this confirms \cite[Conjecture 3.5]{Engel2019ConditionedSystems}.
\begin{mtheorem}[Convergence of finite-time Lyapunov expnonents]
    \label{thm:IntroFTLEconverge}
    Let us assume the same hypotheses as in Theorem~\ref{thm:IntroLyapSpect} such that there exist conditioned Lyapunov exponents $\{\Lambda_i\}_{i=1}^d$. 
    \begin{enumerate}
        \item Then for all $\varepsilon>0$, for $\nu$-almost every $x\in M$,
    \begin{equation*}
    \lim_{t\to\infty}\Pb_x\left[\left|\Lambda_i -\frac{1}{t}\log\delta_i(\Phi_t)\right|>\varepsilon \,\bigg|~ \tau>t \right] = 0. 
    \end{equation*}
    \item If, additionally, for some $p\in (1,\infty]$, we have
\begin{align*} 
    \sup_{t\geq 0}\left\|\frac{1}{t}\log^{+} \left\|\Phi^{\pm1}_t\right\|\right\|_{L^{p}(\Omega \times M, \Pb_\nu(\cdot \mid \tau> t))}<\infty,\label{Qnu}
\end{align*} 
then for $\nu$-almost every $x\in M$
 \begin{equation*}
    \lim_{t\to\infty}\mathbb E_x\left[\left|\Lambda_i -\frac{1}{t}\log\delta_i(\Phi_t)\right| \bigg|~ \tau>t \right] = 0. 
    \end{equation*}
    \end{enumerate}

\end{mtheorem}

{\revision Note that thanks to recent results on QSDs for general diffusions \cite{Benaim2021DegenerateDomain}, our results for SDEs only require hypoellipticity in the sense of the strong H\"{o}rmander condition for the equation of the process $(\varphi_t)_{t\geq 0}$}.

Additionally, we remark that our insights on the random dynamics of the $Q$-process via the invariant ergodic measure $\mathbb Q_{\nu}$  open up the possibility to embark on a programme that resembles the theory of deterministic dynamical systems with holes. In the last two decades, various results have been obtained in~\cite{DemersYoung2006, Bruin2010, DemersWrightYoung2012, DemersTodd2017}, relating the escape rates through holes in a Riemannian manifold, most simply an interval, with the pressure of an invariant ergodic measure on the survival set, taking a similar role as $\mathbb Q_{\nu}$. Here, the pressure is the difference between the metric entropy and the sum of positive Lyapunov exponents with respect to this measure. Upon adding a suitable notion of entropy, we may now be equipped to show similar results for random dynamics with escape. Furthermore, there have been insights on the correspondence between conditionally invariant measures and QSDs \cite{Engel2017LocalDynamics, HomburgZmarrou2007} that may now be extended to their respective relations to survival measures and QEDs. {\revision We mention that, in a similar direction, a quenched approach has been undertaken in \cite{Atnipetal2021, Atnipetal2023ThermoForm, Atnipetal2023EquStates} leading to an ergodic measure on the survival set different to $\mathbb{Q}_{\nu}$.}

The remainder of the paper is structured as follows.  
{\revision Theorem \ref{thm:IntroErgMeas} crucially depends on   Proposition \ref{Q-process} and is proven in Section~\ref{subsec:proof_TheoremA}}; Theorem \ref{thm:IntroLyapSpect}  is contained in the statement of 
Theorem \ref{thm:fullFK} and Theorem~\ref{thm:fullMET}; and Theorem \ref{thm:IntroFTLEconverge} is contained in the statement of {\revision Corollaries \ref{T2.11} and \ref{T2.12} which are consequences of Theorem~\ref{thm:L1Conv}}. In Section \ref{sec:2}, we introduce the setting of this paper and state our main results in more detail.
In Section \ref{sec:3}, we provide results that make the theory of $Q$-process applicable to random dynamical systems with absorption, proving Theorem~\ref{thm:IntroErgMeas} and Theorem~\ref{thm:IntroLyapSpect}.
In Section \ref{sec:3.5}, we link this framework back to finite-time conditioned dynamics proving Theorem~\ref{thm:IntroFTLEconverge}.
In Section \ref{sec:4}, we show how this can be applied to the study of the conditioned dynamics of a large class of stochastic differential equations, significantly generalising the results of \cite{Engel2019ConditionedSystems}.

\section{General setting and main results}
\label{sec:2}
Let $M$ be a $d$-dimensional Riemannian manifold (possibly with boundary) embedded in $\mathbb R^n$. We aim to study random dynamical systems originating inside the domain $M$ and being killed when exiting this region. We denote by $\{\partial\}$ the cemetery state where the flow is absorbed after escape. Moreover, we let ${E_M} := M \sqcup \{\partial\}$ be the topological space generated by the topological basis $\mathcal{T} = \{U; \ U \text{ is open in }M\}\cup \{\partial\}, $ where $\sqcup$ denotes disjoint union. 

Throughout this paper, the time $\mathbb T$ can be taken to be either the semi-group $\mathbb N_0 := \mathbb N \cup \{0\}$ or $\mathbb R_+ := [0,\infty)$. Let  $(\Omega, (\mathcal{F}_t)_{t\in \mathbb{T}}, \mathcal{F}, \Pb, (\theta_t)_{t\in \mathbb{T}})$ be a memoryless noise space (see Appendix \ref{RDS}) where $(\Omega, (\mathcal{F}_t)_{t\in \mathbb{T}}, \mathcal{F})$ fulfils the usual measurability conditions (see \cite[Definition II.67.1]{Rogers2000DiffusionsMartingales}). In this paper, we focus on two different noise spaces, given by 
\begin{equation}\label{eqn:cts}
(\Omega, (\mathcal{F}_t)_{t\geq 0}, \mathcal{F})=\left(\Omega,\left(\sigma((\pi_s)_{0\leq s\leq t})\right)_{t\geq 0}, \sigma((\pi_s)_{s\geq 0})\right) \qquad\mbox{if $\mathbb T = \mathbb R_+$}\,,
\end{equation}
where $\Omega \in \left\{  \mathcal{D}(\R_+,\R^m), \mathcal{C}_0(\R_+,\R^m)\right\}$ or
\begin{equation}\label{eqn:discrete}
(\Omega, (\mathcal{F}_n)_{n\in \mathbb{N}_0}, \mathcal{F})=\left( X^{\mathbb N_0},\left(\sigma((\pi_m)_{0\leq m\leq n})\right)_{n\in\mathbb N_0}, \sigma((\pi_m)_{m\in \mathbb{N}_0})\right) \qquad \mbox{if $\mathbb{T} = \mathbb{N}_0$}\,,
\end{equation}
where $X$ is a Polish space and $\pi_s$ is the canonical processes (see Section \ref{sec:3} for details). These noise spaces are natural for applications to iterated function systems and stochastic differential equations, which we discuss at the end of Section~\ref{secconvprob} and in Section~\ref{Sec2:SDE}. 

Throughout this paper we consider $(\theta,\varphi)$ as a $\mathcal C^1$-random dynamical system on the state space $(E_M,\mathcal B(E_M))$ and with absorption at $\partial$. We further assume that the cocycle $\varphi$ is \emph{perfect} in the sense of Definition \ref{def:rds} (see Appendix \ref{RDS} for details about random dynamical systems).

 \subsection{Absorbed Markov processes}
     
Under the assumption of $ (\Omega, (\mathcal{F}_t)_{t\in \mathbb{T}}, \mathcal{F}, \Pb, (\theta_t)_{t\in \mathbb{T}})$ being a memoryless noise space, $(\theta, \varphi)$ induces a time-homogeneous Markov process $$\varphi = \left(\Omega \times E_M, (\mathcal G_t :=\mathcal{F}_t \otimes \mathcal{B}(E_M))_{t\in \mathbb T}, (\varphi_t)_{t\in\mathbb{T}}, (\mathcal P^t)_{t\in \mathbb T}, (\mathbb{P}_x := \mathbb P\otimes \delta_x)_{x\in E_M}\right)
$$ in the sense of \cite[Definition III.1.1]{Rogers2000DiffusionsMartingales} where
 $\mathcal{P}^t(x,\d y) := \mathbb E_x(\varphi_t \in \d y)$ for every $x\in E_M$, i.e.~
\begin{itemize}
    \item[(i)] $\left(\Omega \times E_M, \mathcal{G}_t,\mathcal G \right)$ is a filtered space, where $\mathcal G := \sigma \left( (\mathcal G_t)_{t\in \mathbb T}\right)$. From equations \eqref{eqn:cts} and \eqref{eqn:discrete} we have that $\mathcal G = \mathcal F \otimes \mathcal B(E_M)$;
    \item[(ii)] $\varphi_t$ is an $\mathcal G_t$-adapted process with state space $E_M$;
    \item[(iii)] $\mathcal P^t$ a time-homogeneous transition probability function of the process $\varphi_t$ satisfying the usual measurability assumptions and the Chapman--Kolmogorov equation;
    \item[(iv)] $(\mathbb P_x)_{x\in{E_M}}$ is a family of probability measures satisfying $\mathbb P_x[\varphi_0=x] = 1$ for every $x\in E_M$; and
    \item[(v)] for all $t,s\in \mathbb T$ and  every bounded measurable function $f$ on $M$ 
 $$\mathbb E_x\left[f\circ \varphi_{t+s}\mid \mathcal G_t\right] = ({\mathcal P}^s f)(\varphi_t)  \qquad { \ \mathbb P_x\text{-almost surely}}.$$
\end{itemize} 
For a proof,  refer to \cite[Chapter 2.5]{Newman2015SynchronisationSystems}.

Since $\varphi_t$ is absorbed at $\partial,$ 
we can define the stopping time
    $$\tau(\cdot , x) = \inf\left\{t\geq 0 : \varphi_t(\cdot, x) = \partial\right\}.$$

Below, we introduce some notations used throughout the present paper. 

\begin{notation}
Given a measure $\mu$ on $M$, we denote
$\mathbb P_\mu (\cdot) := \int_{M} \mathbb P_x (\cdot) \mu (\d x).$

 We consider the set $\mathcal F_b(M)$ as the set of bounded Borel measurable functions on $M$. Given $f\in \mathcal F_b(M)$ and $g\in \mathcal F_b(\Omega\times M)$, by abuse of notation we write
\begin{equation} \label{eq:semigroup}
   \mathcal P^t(f)(x):= \mathcal P^t\left(\mathbbm 1_M   f \right)(x) = \int_M f(y) \mathcal P^t(x, \mathrm{d}y), 
\end{equation} 
$$\mathbb E_x \left[g\right] := \mathbb E_x \left[ \mathbbm 1_M  g\right], \ {\   \text{for all}} \ x\in M, $$
and
$$ f\circ \varphi_t := (\mathbbm 1_M \circ \varphi_t)\cdot (f\circ \varphi_t)$$

We denote by $\mathcal C^0(M)$ the space of continuous functions $f:M\to M$, and by $\mathcal M(M)$ the set of Borel signed measures on $M$. 

 An essential tool to the study of random dynamical systems is invariant measures
 which correspond with stationary measures of the associated Markov process.
 However, in the context of absorbed dynamics, such measures do not exist due to the exponential loss of mass of $\varphi$ on $M$. These measures are replaced by so-called \emph{quasi-stationary distributions} $\mu$ (QSD).
\end{notation}

\begin{definition}[Quasi-stationary distribution] A probability measure $\mu$ on $(M, \mathcal{B}(M))$ is said to be a \emph{quasi-stationary distribution} for the random dynamical system $(\theta,\varphi)$ if for all $A \in \mathcal{B}(M)$
$$\Pb_{\mu}(\varphi_t \in A \, |\, \tau>t)= \mu(A) \qquad \text{for all }t\in\mathbb T.$$
\end{definition}
Note that, since $\varphi_t$ is absorbed at $\partial$, we have
$$\Pb_{\mu}(\varphi_t \in A \, |\, \tau>t) = \frac{\mathbb P_\mu(\varphi_t \in A)}{\Pb_\mu (\tau >t)} =  \frac{\int_{M} \mathcal P^t(x,A) \mu(\d x)}{\int_{M} \mathcal P^t(x,M) \mu(\d x)}\qquad \text{for all } t\in\mathbb T.$$

Furthermore, if the absorbed dynamics evolve under the statistics of a unique quasi-stationary distribution, in contrast, the history of the surviving trajectories at time $T>0$ do not, in general, follow the quasi-stationary statistics. Instead, the asymptotic distribution of the history of surviving trajectories is given by the so-called \emph{quasi-ergodic distribution} $\nu$ (QED).

\begin{definition}[Quasi-ergodic distribution] \label{Def:QEM}
A probability measure $\nu$ on $(M, \mathcal{B}(M))$ is said to be a \emph{quasi-ergodic distribution} for  the random dynamical system $(\theta,\varphi)$ if for all $f : M \to \R$ bounded and $\mathcal{B}(M)$-measurable
$$\begin{cases}\displaystyle
\lim_{t\to\infty}\mathbb{E}_x \left[\frac{1}{t}\int_0^t f(\varphi_s) \mathrm{d} s \Bigg|~ \tau>t\right] = \int_M f \ \mathrm{d} \nu \qquad \mbox{ for all } x\in M,& \ \text{if }\mathbb T = \mathbb R_+,\\
\displaystyle \lim_{n\to\infty}\mathbb{E}_x \left[\frac{1}{n}\sum_{i=0}^{n-1} f(\varphi_i)  \Bigg|~ \tau>n\right] = \int_M f\ \mathrm{d} \nu \qquad  \mbox{ for all } x\in M,& \ \text{if }\mathbb T =\mathbb N_0.
\end{cases}$$

\end{definition}

{\revision
Recall from the Introduction that for this paper we impose Hypothesis \ref{(H)}, ensuring the existence and uniqueness of a QSD and QED.}
In particular, we require the RDS with absorption $\varphi$ to have pointwise exponential convergence towards the QSD in the total variation norm. 

For criteria ensuring this hypothesis, see \cite{Benaim2021DegenerateDomain, Castro2021ExistenceApproach, Champagnat2018CriteriaDiffusions, Champagnat2016ExponentialQ-process}. Some properties of the conditioned process induced by these conditions are given by Proposition \ref{prop:A}, {\revision in particular the relation
\begin{equation}
\label{eq:QSD_QED}
\nu(\d x) = \eta(x) \mu(\d x),
\end{equation}
where $\eta$ is the bounded function in ${\mathrm{(H3)}}$.
}
Notably, our setting ensures the existence of the  $Q$-process \cite{Champagnat2016ExponentialQ-process} in the strongest possible sense, a key element to the proof of the multiplicative ergodic theorem in the conditioned setting.

The literature on absorbed Markov processes generally assumes the conditions given by  \cite{Champagnat2016ExponentialQ-process} which imply exponential convergence of \eqref{lim1} uniformly on $x$. While these conditions are well suited for the study of stochastic differential equations with escape, the uniform convergence with respect to $x$ of $\eqref{lim1}$ turns out to be too restrictive for discrete-time systems with escape, specifically with bounded noise (see \cite{Castro2021ExistenceApproach, castro_goverse_lamb_rasmussen_2023}).

{\revision 
\begin{remark} We point out some important perspectives on $\mu$, $\eta$ and $\nu$:
\begin{enumerate}
    \item[(i)] We can obtain $\eta$ as dominant eigenfunction of the sub-Markovian semigroup $\mathcal P^t$ \eqref{eq:semigroup}, i.e.~$\mathcal P^t \eta = e^{- \beta t} \eta $, where $\beta$ is the escape rate with respect to the QSD $\mu$. Analogously, we have that $\mu$ is an eigenmeasure of the adjoint semigroup for the same eigenvalues, i.e.~$(\mathcal P^t)^* \mu = e^{- \beta t} \mu $ (see e.g.~\cite{Champagnat2016ExponentialQ-process}).
    Note that in the case with no killing, we have $\beta =0$ such that $\eta$ is simply a constant and $\mu$ is a solution of the stationary forward Kolmogorov problem (if it exists).

    \item[(ii)] When there is loss of mass through the cemetery state, the function $\eta$ is typically non-constant and therefore expresses a discrepancy between the QSD $\mu$ and the QED $\nu$.
    For example, consider the most simple SDE $\d X_t = \d W_t$ on $[0,\pi]$, corresponding with the killed heat semigroup $\mathcal P^t$ and generator $\Delta$ with Dirichlet boundary conditions.
    We have $\mu(\d x) = \sin(x)/2 \d x$ and $\eta(x) = (4/\pi)\sin(x)$ such that $\nu(\d x) = (2/\pi) \sin^2x \d x$. Clearly, compared to $\mu$, the QED $\nu$ has stronger concentration around the center of the interval $[0,\pi]$ accounting for the fact that the Birkhoff sums collect the whole history of trajectories not being killed.
    

    \item[(iii)] It turns out (cf.~\cite[Proposition 6.4.2]{Engel2017LocalDynamics} or \cite[Lemma 5.2]{HomburgZmarrou2007}) that the measure $\mathbb P \otimes \mu$ is a conditionally invariant measure of the skew product flow $\Theta_t := (\theta_t, \varphi_t)$ on $M\times \Omega$ (see Theorem~\ref{thm:IntroErgMeas} and below for more details) if and only if $\mu$ is a QSD for $\varphi_t$. 
    
    Apart from this observation on the skew product flow setting, for general deterministic systems with holes, conditionally invariant measures describe the stationary statistics in analogy to QSD (see \cite[Chapter 8]{Collet2013Quasi-StationaryDistributions}).

    One may also consider a similar analogy between a QED $\nu$ and an invariant measure defined on the infinite-time survival
    set of a deterministic system with a hole. Observe that for bounded and measurable $f$ we have
    $$\lim_{t\to\infty}e^{\beta t}\mathbb{E}_{\mu}\left[f(\varphi_0)\mathbbm{1}_{\{\tau>t\}}\right] = \lim_{t\to\infty}\int_{M}f(x)e^{\beta t}\mathbb{P}_x(\tau>t)\mu(\d x) = \int_M f(x)\eta(x)\mu(\d x) = \nu(f)$$
    where we have used \eqref{eqn:etalim}. The left-hand side is analogous to the characterisation of the above-mentioned invariant measure (see for instance \cite[Theorem 2.16]{Bruin2010}). The discrepancy between $\mu$ and $\nu$ is thus similar to the one observed for deterministic systems with holes.
\end{enumerate}
    
\end{remark}
}

%

%

%
%

%
%
%
%
%
%
%
%
 
%
%
Under  Hypothesis \ref{(H)}, we may further assume the existence of the $Q$-process shown by \cite{Champagnat2016ExponentialQ-process, Champagnat2017UniformQ-process},  the process $(\varphi_t)_{t\in \mathbb{T}}$ conditioned on asymptotic survival.

\begin{definition}[$Q$-process] 
A family of probability measures $(\mathbb Q_x)_{x \in M}$ on $\left(\Omega \times M,\mathcal G \right)$ is called $Q$-process, if 
\begin{enumerate}
    \item for each $x \in M$ and every $s\geq 0$ and $A \in \mathcal G_s$ we have
    \begin{align}
        \mathbb Q_x (A) := \lim_{t\to\infty} \mathbb P_x (A\mid \tau>t) \label{Qlimit:2.6}\,,
    \end{align} 
    
    \item the tuple $$\left(\Omega\times M, \left(\mathcal{G}_t\right)_{t\in\mathbb T}, \left(\varphi_t\right)_{t\in\mathbb T}, \left(\mathcal Q^t\right)_{t\in\mathbb T}, (\mathbb Q_x)_{x \in M}\right)$$ is a Markov process, where we define
    $$\mathcal Q^t(x,B) := \mathbb Q_x[\varphi_t \in B]\ \ \text{for all }t\in\mathbb T,\ x\in M \text{ and }B\in\mathcal B(M).$$
\end{enumerate}
\label{def:Qproc}
\end{definition}
Note that, by definition, a $Q$-process is unique. 

    In previous descriptions of the $Q$-process such as \cite{Champagnat2016ExponentialQ-process,Champagnat2017UniformQ-process}, the probability measures $(\mathbb Q_x)_{x\in M}$ were only defined on $\bigcup_{t \geq 0} \mathcal{G}_t$. However, for our application to random dynamical systems, it is important for the measures $(\mathbb Q_x)_{x\in M}$ to be extended to  $\mathcal G$. 
    
    To understand why one cannot expect the limit \eqref{Qlimit:2.6} to hold for all $A \in \mathcal G$, consider the set $\left\{  \tau  = \infty\right\}\in \mathcal G.$
    Under Hypothesis \ref{(H)}, we have $\mathbb P_x(\tau = \infty) = 0$
    for all $x \in M$ and thus also $\mathbb P_x(\tau = \infty\mid \tau >t) = 0, \text{ for all $x\in M$ and $t>0$.}$
    On the other hand,
    $$\mathbb Q_x (\{\tau > s\}) = \lim_{t\to \infty}\mathbb P_x(\tau > s |\tau > t) = 1, \text{ for all $x\in M$ and $s>0$.}$$
    and thus $$\mathbb{Q}_x(\tau = \infty) = 1 \neq 0 = \lim_{t \to \infty}\mathbb P_x (\tau = \infty\mid \tau >t).$$

\begin{proposition}[Existence of the $Q$-process]\label{Q-process} Under Hypothesis \ref{(H)}, there exists a $Q$-process  $(\mathbb{Q}_x)_{x \in M}$ with transition kernels given by
$$\mathcal Q^t(x,\mathrm{d}y) = e^{\beta  t}\frac{\eta(y)}{\eta(x)}\mathcal{P}^t(x, \mathrm{d}y).$$

Furthermore, the measure $\nu$ is the unique stationary measure of the Markov process $$\big(\Omega\times M, \left(\mathcal{G}_t\right)_{t\in\mathbb T}, \left(\varphi_t\right)_{t\in\mathbb T}, \left(\mathcal Q^t\right)_{t\in\mathbb T}, (\mathbb Q_x)_{x \in M}\big)$$ and we have
$$\lim_{t\to\infty}\left\|\mathbb{Q}_x(\varphi_t \in \cdot) - \nu\right\|_{TV}=0.$$
\end{proposition}

Proposition \ref{Q-process} is proved in Section \ref{subsec:Q-process}.

\subsection{{\revision The random dynamical system on the survival set and Theorem~\ref{thm:IntroErgMeas}}}
Recall that the skew product $(\Theta_t)_{t\in \mathbb{T}}$ of $(\theta,\varphi)$ defined as
\begin{align*}
    \Theta_t : \Omega \times E_M &\to \Omega \times E_M\\
    (\omega, x) &\mapsto \Theta_t(\omega, x) := \left(\theta_t \omega, \varphi(t, \omega, x)\right)
\end{align*}
is a family of measurable mappings generating a semi-flow, i.e.~a measurable dynamical system. 

The proof of the existence of conditioned Lyapunov exponents relies on finding a suitable ergodic probability measure for $\Theta_t$ giving full measure to paths never to be absorbed, i.e.~to
\begin{equation}
\Xi := \left\{(\omega, x) \mid \tau(\omega, x) = \infty\right\} = \bigcap_{t\in \mathbb T} \left\{(\omega,x) \in \Omega\times M\, \mid \tau(\omega,x)> t\right\}\label{eqn:Xi}.
\end{equation}

The most appropriate choice for such a measure turns out to be
$$\mathbb Q_\nu (\cdot) := \int_{M} \mathbb Q_x (\cdot) \nu (\d x).$$
Observe that the measure $\mathbb Q_{\nu}$ on $(\Omega\times M, \mathcal F\otimes \mathcal B(M) )$ satisfies
\begin{align}
    \mathbb{Q}_{\nu}(\Xi) = 1\label{Qv}
\end{align}
and $\Pb_x(\Xi) = 0,$ for every $x\in M.$ Given a function $f\in L^1(M\times \Omega, \mathbb Q_\nu),$ we denote
$$\mathbb E_\nu^\mathbb Q [f] = \int_{\Omega \times M } f(\omega,x) \ \mathbb Q_{\nu} (\d \omega, \d x).$$

Note here that we impose the assumption of a perfect cocycle as for all $x \in M$. However, this can be loosened: if the cocycle is not perfect, even under the $Q$-measures, the cocycle property holds almost surely (see Remark \ref{rmk:modification} below).

{\revision
As presented in the Introduction, the first central statement and a crucial insight of this paper is given by Theorem~\ref{thm:IntroErgMeas} which states that the measure $\mathbb Q_\nu$ is ergodic and, hence, fulfils an essential condition for the proof of the existence of conditioned Lyapunov exponents. The proof can be found in Section~\ref{subsec:proof_TheoremA}.
}

\subsection{Conditioned Lyapunov exponents}\label{CondLyapunovExp}

Since $M$ is a manifold embedded in $\mathbb R^n$, we can consider $$TM=\left\{(x,v) \in M\times \R^n, x\in M\ \text{and }v\in T_x M\right\},$$
as the tangent bundle of $M,$ where $T_x M$ denotes the tangent space of $M$ at $x$ (for a complete description of the tangent bundle and its properties see \cite[page 65]{Lee2012IntroductionManifolds}). From 
Theorem~\ref{thm:IntroErgMeas} and equation \eqref{Qv} we have that $\mathbb Q_\nu$ is an ergodic measure to the dynamical system $(\Theta_t)_{t\in\mathbb T}$ and $\mathbb{Q}_{\nu}(\Xi) = \mathbb Q_\nu [\tau = \infty] = 1.$

Observe that for $\mathbb Q_\nu$-a.e.~$(x,\omega)\in \Omega\times M$ the linear map
\begin{align*}
    \Phi_t(\omega,x):T_{x}M&\to  T_{\varphi_t(\omega,x)}M\\
    v&\mapsto \mathrm{D} \varphi_t(\omega,x)v
\end{align*} 
is well-defined where $\mathrm{D}$ denotes the derivative w.r.t.~$x$. 
Throughout this paper, $\Phi_t$ might refer to the map above or be considered acting on the state space $TM$. 

Moreover, since $\varphi_{t+s}(\omega,x) = \varphi_t(\theta_s\omega, \varphi_s(\omega,x)),$ for every $t,s\in\mathbb T,$ by differentiation in $x$ and using the chain rule we obtain 
$$ \Phi_{t+s}(\omega,x) = \Phi_t(\Theta_{s}(\omega,x))\circ \Phi_s(t,\omega), \quad \text{ for } \mathbb Q_\nu\text{-almost every }(\omega,x)\in \Omega\times M,$$
i.e.~$\Phi_t$ forms a cocycle over the dynamical system $(\Theta_t)_{t\in\mathbb T}.$

For $\mathbb{Q}_{\nu}$-almost every $(\omega, x) \in \Omega \times X$, we wish to show the  convergence of the following limits $\Lambda_i(\omega, x)$ defined by
$$\Lambda_i(\omega, x) = \lim_{t\to \infty}\frac{1}{t}\log\delta_i(\Phi_t(\omega, x)) \qquad \mbox{for all } i \in\{ 1, \ldots, d\}\,,$$
where $\delta_i(\Phi_t)$ denotes the $i\textsuperscript{th}$ singular value of $\Phi_t$, i.e.~the %
square root of the $i\textsuperscript{th}$ eigenvalue of $\Phi_t^{*}\Phi_t$, when $\Phi_t$ is seen as an $\R^d$-endomorphism. For $(\omega, x, v) \in \Omega \times TM$, one may 
define the finite-time Lyapunov exponents $\lambda_v(t, \omega, x)$
\begin{align*}
    \lambda_v(t,\omega, x) = \frac{1}{t}\log\frac{\left\|\Phi_t(\omega, x)v\right\|}{\|v\|} 
\end{align*}
where $\|\cdot\|$ is induced by the Riemannian metric on $M$ and their limit superiors $ \lambda_v(\omega, x)$, the 
\emph{characteristic Lyapunov exponents}
\begin{align*}\lambda_v(\omega, x) = 
\limsup_{t\to\infty}
\frac{1}{t} \log\frac{\left\|\Phi_t(\omega, x)v\right\|}{\|v\|}.
\end{align*}
We observe that, in our setting,  $\lambda_v(\omega, x)$ exists as an actual limit and takes one of the values of $\Lambda_i(\omega,x)$.
In fact, it is directly related to the ergodicity of $(\Theta_t)_{t\in \mathbb{T}}$ with respect to $\mathbb{Q}_{\nu}$ and an application of the Furstenberg--Kesten theorem, that the RDS associated with the $Q$-process
has a spectrum of Lyapunov exponents, which do not depend on initial conditions $(\omega,x).$

\begin{theorem}[Spectrum of Lyapunov exponents] \label{GrasmannianMeasure} Let $(\Phi_t)_{t\in \mathbb T}$ be as above, i.e.~an RDS over the dynamical system $\left(\Omega \times M, \mathcal{F}\otimes\mathcal{B}(M), (\Theta_t)_{t\in \mathbb{T}}\right)$ %
with ergodic invariant measure
$\mathbb{Q}_{\nu}$. Assume further that
\begin{equation}
\mathbb{E}^{\mathbb{Q}}_{\nu}\left[\sup_{0\leq t\leq 1}\log^{+}\|\Phi_t\|\right]<\infty\label{eqn:IC+}.
\end{equation}
\begin{enumerate}
    \item Then there exists a $\Theta$-forward-invariant set $\Delta \in \mathcal{F}\otimes \mathcal{B}(M)$ of $\mathbb{Q}_{\nu}$-full measure such $\Delta \subset \Xi$ (see \eqref{eqn:Xi}) and constant Lyapunov exponents ${\Lambda_1\geq \ldots\geq \Lambda_d\geq -\infty}$ such that for all ${(\omega, x) \in \Delta}$
    \begin{equation}
    \Lambda_i = \lim_{t\to \infty}\frac{1}{t}\log\delta_i(\Phi_t(\omega, x)) \qquad \mbox{for all } i \in\{ 1, \ldots, d\} \label{eqn:LyapConv}
    \end{equation}
    \item If in addition, $\Phi_t(\omega, x)$ is invertible for all $(t, \omega, x) \in \mathbb{T}\times \Omega \times M$ and
    \begin{equation}\mathbb{E}^{\mathbb{Q}}_{\nu}\left[\sup_{0\leq t\leq 1}\log^{+}\|\Phi^{-1}_t\|\right]<\infty,\label{eqn:IC-}\end{equation}
    then the Lyapunov exponents are finite and the convergence \eqref{eqn:LyapConv} holds in $L^1(\mathbb{Q}_{\nu}).$
    \item Let $\lambda_i, d_i$ respectively denote the distinct Lyapunov exponents and their multiplicities; let $p$ be the number of distinct Lyapunov exponents. Then we define the Lyapunov spectrum
    $$\mathcal{S}(\theta, \varphi) = \left\{(\lambda_i, d_i) : i = 1, \ldots, p\right\}.$$

\end{enumerate}
\label{thm:fullFK}
\end{theorem}

Furthermore, the setting of the $Q$-process also yields the existence of Oseledets flags.

\begin{theorem}[Multiplicative ergodic theorem] Let $(\Theta, \Phi)$, $\Delta$, $\mathcal{S}(\theta, \varphi)$ be as above. Then for all~$(\omega, x) \in \Delta$, the following statements hold:
\begin{enumerate}
    \item the random matrix limit $\Psi(\omega, x) := \lim_{t\to\infty}\left(\Phi_t(\omega, x)^{*}\Phi_t(\omega, x)\right)^{1/2t}$ exists and has eigenvalues ${e^{\lambda_1}>\ldots > e^{\lambda_p}}.$
    \item Let $E_1(\omega, x),\dots, E_p(\omega, x) \subset T_x M$ denote the corresponding random eigenspaces of $\Psi(\omega, x)$ with ${\dim E_i = d_i}$ and define
    $$U_i(x,\omega) = \bigoplus_{k=i}^p E_k(\omega, x).$$
    Then $\Phi_t(\omega, x)U_i(\omega, x) = U_i(\Theta_t(\omega, x))$ and the $U_i$'s form a random filtration of $T_xM$:
    $$\{0\}\subset U_p(\omega, x) \subset U_{p-1}(\omega, x)\subset\ldots \subset U_2(\omega, x) \subset U_1(\omega, x) = T_x M.$$
    \item Furthermore for all $v\in T_x M$, the finite-time Lyapunov exponents converge and   
    $$\lambda_v(\omega, x) = \lim_{t\to\infty}\frac{1}{t}\log\|\Phi_t(\omega, x) v\| = \lambda_i \Longleftrightarrow v \in U_i(\omega, x) \backslash U_{i+1}(\omega, x).$$
\end{enumerate}
\label{thm:fullMET}
\end{theorem}

These two theorems are obtained directly from Theorem~\ref{thm:IntroErgMeas} in combination with the classical theory of random dynamical systems (see Section \ref{QMET} for more details). This shows that the $Q$-process setting is well-suited to studying conditioned dynamics. However, we wish to ensure in Section \ref{sec:3.5} that in the particular context of absorbed diffusion processes, this corresponds precisely to the framework introduced by Engel et al.~\cite{Engel2019ConditionedSystems}.

In \cite{Engel2019ConditionedSystems}, the existence of the top Lyapunov exponent  is proved by introducing an extended process on the unit tangent bundle. This can be generalised as in \cite{Baxendale1986TheDiffeomorphisms} for the full spectrum of Lyapunov exponents by introducing a process on the Grassmannian bundle $\mathrm{Gr}_k(M)$, whose fibers $\mathrm{Gr}(T_xM)$ are the manifolds consisting of subspaces of the tangent spaces $T_x M$ (see Section \ref{METQP}). Defining the space of the alternating $k$-multivectors
$$\Bigwedge^k_0 T_x M = \left\{v_1 \wedge \cdots \wedge v_k \mid v_1,\ldots, v_k \in T_x M \right\}$$
generating the vector space $\Bigwedge^k T_x M$. One can equivalently identify $\mathrm{Gr}_{k} (T_x M)$ as the set $\mathbf{P}(\Bigwedge^k_0 T_x M)$, which is a $k(d-k)$-dimensional submanifold of the projective space $\mathbf{P}(\Bigwedge^k T_x M)$. Furthermore, let us define the vector space homomorphism
\begin{align*}
    \Bigwedge^k \Phi_t(\omega, x) : \Bigwedge^k T_x M &\to \Bigwedge^k T_{\varphi_t(\omega, x)}M
\end{align*}
defined on $\Bigwedge^k_0 T_x M$ by
$$\Bigwedge^k \Phi_t \left(v_1 \wedge \cdots \wedge v_k\right) := \Phi_t(\omega, x) v_1 \wedge \cdots \wedge \Phi_t(\omega, x) v_k. $$
See for instance \cite{Crauel1990LyapunovGrassmannians} for a concise introduction of exterior powers in the context of Lyapunov exponents. This allows us to state the following proposition.

\begin{proposition}\label{thm:L1QConv}
    Assume that $(\Theta, \Phi)$ satisfies integrability conditions \eqref{eqn:IC+} and \eqref{eqn:IC-} and for $k\leq d$, let $ \rho^k$ be the Borel measure on $\mathrm{Gr}_k(M)$ defined as
    $$ \rho^k(  \d  x \times \d v) = \sigma_x^k\left( \d v \cap \mathrm{Gr}_k(T_x M)\right)\nu(\d x),$$ 
    where  $\sigma^k_x$ is the Haar measure on $\mathrm{Gr}_k(T_x M)$ defined in \cite[p. 325]{Baxendale1986TheDiffeomorphisms}. Then there exists a set ${\widetilde{G} \subset \mathrm{Gr}_k(M),}$ such that $\rho^k(\widetilde{G}) = 1$ and
    \begin{equation*}
    \lim_{t\to\infty}\mathbb{E}^{\mathbb{Q}}_{\nu}\left[ \left|\lambda^{(k)} - \frac{1}{t}\log\left\|\Bigwedge^k \Phi_tv\right\|\right| \right] = 0\qquad \text{for all }(x,v)\in \widetilde{G},
    \end{equation*}
    where $\lambda^{(k)} = \Lambda_1 + \cdots +\Lambda_k$ and $\left\{\Lambda_i\right\}_{i =1}^d$ are given by Theorem \ref{thm:fullFK}.
\end{proposition}

Note that the integrability conditions \eqref{eqn:IC+} and \eqref{eqn:IC-} of the multiplicative ergodic theorem can be formulated in terms of the conditioned process. Indeed from Proposition \ref{marathon} below, with the QSD $\mu$, we have
\begin{align} \label{ineq:QSDcond}
    \mathbb{E}_{\nu}^{\mathbb Q} \left[\sup_{0\leq t\leq 1}\log^+ \|\Phi^{\pm1}_t \| \right] &=  \int_M \mathbb{E}_{x}^{\mathbb Q} \left[\sup_{0\leq t\leq 1}\log^+ \|\Phi^{\pm1}_t \| \right]  \nu(\d x) \nonumber\\
    &= e^{\beta }\int_M \mathbb  E_x \left[\eta(\varphi_1)\sup_{0\leq t\leq 1}\log^+ \|\Phi^{\pm1}_t \| \right]  \mu(\d x) \nonumber\\
    &\leq  e^{\beta }\|\eta\|_{\infty} \mathbb E_\mu \left[\sup_{0\leq t\leq 1}\log^+ \|\Phi^{\pm1}_t \| \mathbbm 1_{\{\tau >1\} } \right].
\end{align}

\subsection{Convergence in conditional probability}\label{secconvprob}

\begingroup
\renewcommand\thehypothesis{(C)}

Under the same setting and mild conditions, we prove the convergence of the finite-time Lyapunov exponents towards the $Q$-process Lyapunov exponents in conditional probability. To this end, we exhibit the following result in the more general settings of Markov processes with $Q$-processes: under suitable conditions, we show that any convergence $\Gamma_t\to\Gamma^{\star}$ in mean or in probability under $\mathbb{Q}_x$ also holds respectively in conditional mean or probability under $\Pb_x$. For Theorem \ref{thm:L1Conv} below, although we keep the notation $(\varphi_t)_{t\in\mathbb{T}}$, it does not necessarily denote a cocycle but any absorbed Markov process $(\widetilde{\Omega},(\mathcal{G}_t)_{t\in \mathbb{T}},(\varphi_t)_{t\in\mathbb{T}}, (\mathcal{P}^t)_{t\in\mathbb{T}}, (\mathbb{P}_x)_{x\in M\sqcup\left\{\partial\right\}})$ for which there is a corresponding $Q$-process $(\widetilde{\Omega},(\mathcal{G}_t)_{t\in \mathbb{T}},(\varphi_t)_{t\in\mathbb{T}}, (\mathcal{Q}^t)_{t\in\mathbb{T}}, (\mathbb{Q}_x)_{x\in M})$ under Hypothesis \ref{(H)}.

{\revision The following theorem states the other main insight of this paper yielding the two subsequent corollaries that are summarised in Theorem C. It allows the identification of the Lyapunov exponents obtained from the $Q$-process description with the limits of conditioned finite-time Lyapunov exponents:}

\begin{theorem}\label{thm:L1Conv}
    Let $(\varphi_t)_{t\in \mathbb{T}}$ be a Markov process satisfying Hypothesis \ref{(H)}. Let $x\in M$ and $\left(\Gamma_t\right)_{t\in \mathbb{T}}$ be a collection of $\mathcal G_t$-measurable random variables.
\begin{itemize}
    \item[(i)] Suppose that $(\Gamma_t)_{t\in \mathbb T}$ convergences in probability to some $\Gamma^{\star} \in \R$ under $\mathbb{Q}_x$, i.e.~for all $\varepsilon > 0$,
    \begin{equation*}
        \lim_{t\to \infty}\mathbb{Q}_x\left[\left|\Gamma_t - \Gamma^{\star}\right|>\varepsilon\right] = 0.
    \end{equation*}
    Then this convergence holds in $\Pb_x$-conditional probability, i.e.~for all $\varepsilon>0$,
    \begin{equation*}
    \lim_{t\to \infty}\mathbb{P}_x\left[\left|\Gamma_t - \Gamma^{\star}\right|>\varepsilon\mid \tau>t\right] = 0.
    \end{equation*}
  
    \item[(ii)] If in addition, there exists $p \in (1, \infty]$ such that
    \begin{equation}
    \lim_{t\to \infty} \mathbb{E}_{x}^{\mathbb Q}[|\Gamma_t - \Gamma^{\star}|] =0 \qquad \textrm{and}\qquad \sup_{t\geq 0} \| \Gamma_t\|_{L^p(M\times \Omega, \mathbb P_x(\cdot \mid \tau>t) )}<\infty \label{normal}
    \end{equation}
    then 
    \begin{equation*}
        \lim_{t\to \infty} \mathbb E_x [|\Gamma_t - \Gamma^{\star}|\mid \tau >t] = 0.
    \end{equation*}
\end{itemize}
\end{theorem}

Using these insights,
we obtain the following convergence theorems for 
finite-time Lyapunov exponents 
under conditional probabilities. 

\begin{corollary}\label{T2.11}
Assume that $(\Theta, \Phi)$ fulfils the integrability condition \eqref{eqn:IC+} so that the multiplicative ergodic theorem holds. Let $k\leq d$ be such that $$\lambda^{(k)}=\Lambda_1+\cdots +\Lambda_k >- \infty$$
where $\Lambda_1,\ldots, \Lambda_k$ are the $k$ first Lyapunov exponents given by Theorem \ref{thm:fullFK}.
    Then for every $\varepsilon> 0$ and $\rho^k$-almost every $(x,v) \in \mathrm{Gr}_k(M),$
\begin{align}
    \lim_{t\to\infty } \mathbb P_x \left[ \left\{ \left|\lambda^{(k)} - \frac{1}{t}\log \|\Bigwedge^k \Phi_t v\|\right| >\varepsilon \right\}\Bigg|~  \tau > t \right]  = 0.\label{eqn:CondProbV}
\end{align}

Similarly, for all $\varepsilon>0$ and $\nu$-almost every $x \in M$

\begin{equation}\lim_{t\to\infty}\mathbb P_x\left[ \left\{\left|\lambda^{(k)} - \frac{1}{t} \log\left\|\Bigwedge^k \Phi_t \right\|\right|> \varepsilon \right\}\bigg|~ \tau>t\right] = 0.\label{eqn:CondProbNorm}\end{equation}

\end{corollary}
Note that, due to inequality~\eqref{ineq:QSDcond}, the crucial condition~\eqref{eqn:IC+} (and similarly \eqref{eqn:IC-}) follow readily from
$$ \mathbb E_\mu \left[\sup_{0\leq t\leq 1}\log^+ \|\Phi^{\pm1}_t \| \mathbbm 1_{\{\tau >1\} } \right] < \infty,$$
which is easier to verify in practice by explicit knowledge of the QSD $\mu$.
\begin{remark}
    Observe that a slightly different version of Corollary \ref{T2.11} is the following: For all $v \in T_x M$ such that $\lambda_v :=\lambda(\cdot, \cdot, v) : (\omega, x) \mapsto \lambda(\omega, x, v)$ is constant $\mathbb{Q}_{\nu}$-almost surely, there exists $k\leq p\leq d$ such that $
    \lambda(\cdot, \cdot, v)= \lambda_k$ holds $ \mathbb{Q}_{\nu}$-almost surely. Furthermore,
    \begin{align}
    \lim_{t\to\infty } \mathbb P_x \left[\left|\lambda_k - \frac{1}{t}\log \left\|\Phi_t v\right\| \right|>\varepsilon \bigg|~    \tau > t \right]  = 0. \label{C1}
\end{align}
This is useful in cases where some of the Oseledet's flags are not random or degenerate (see Example~\ref{ex:uncoupled} below). As a matter of fact, it is believed that the Oseledet's spaces are either constant or that their distribution is non-degenerate, at least for a large class of stochastic differential equations. This would immediately imply that the conditioned characteristic Lyapunov exponents $\lambda_v$ are constant $\mathbb{Q}_{\nu}$-almost surely for all $v\in  T_x M$. However, we were not able to find such known general result that would most likely rely on the use of Malliavin calculus in the spirit of \cite{Imkeller1998TheEquations}.
\end{remark}

Under stronger assumptions, the convergences \eqref{eqn:CondProbV} and \eqref{eqn:CondProbNorm} can be strengthened in conditional expectation.
\begin{corollary}\label{T2.12}
Assume now that $(\Theta, \Phi)$ is an invertible linear cocycle fulfilling the integrability conditions \eqref{eqn:IC+} and $\eqref{eqn:IC-}$ and let $\rho^k$ be as above. Assume further that for some $p\in (1,\infty]$
\begin{align} 
    \alpha^{\pm} = \sup_{t\geq 0}\left\|\frac{1}{t}\log^{+} \left\|\Phi^{\pm1}_t\right\|\right\|_{L^{p}(\Omega \times M, \Pb_\nu(\cdot \mid  \tau> t))}<\infty.\label{Qnu}
\end{align} 

Then $\lambda^{(d)}>-\infty$ and for all $k\leq d$, for $\rho^k$-almost every $(x,v) \in \mathrm{Gr}_k(M),$
\begin{align}
    \lim_{t\to\infty } \mathbb E_x \left[\left|\lambda^{(k)} - \frac{1}{t}\log \left\|\Bigwedge^k\Phi_t v\right\| \right|\bigg|~     \tau > t \right]  = 0\label{C4}
\end{align}

and for $\nu$-almost every $x\in M$
\begin{equation}
\lim_{t\to\infty}\mathbb E_x\left[\left|\lambda^{(k)} - \frac{1}{t}\log \left\|\Bigwedge^k\Phi_t\right\| \right|\bigg|~   \tau>t\right] = 0.\label{C5}
\end{equation}
\end{corollary}

Theorem \ref{thm:L1Conv} and {\revision Corollaries \ref{T2.11} and \ref{T2.12}} are proved in Section \ref{sec:3.5}.

\begin{example}[Iterated function systems]\label{ex:IFS} Let $K$ be a compact metric space and $\Pi$ be a Borel probability measure on $K$; then they generate a memoryless noise space
$$\left( K^{\mathbb{N}_0}, \mathcal B (K)^{\otimes \mathbb N_0}, (\mathcal{F}_n)_{n\in \mathbb{N}_0}, (\theta_n)_{n\in \mathbb{N}_0}, \Pi^{\otimes \mathbb N_0} \right)$$
of the form \eqref{l2}, where
\begin{align*}
    \theta := \theta_1 : K^{\mathbb{N}_0} &\rightarrow K^{\mathbb{N}_0}\\
    \omega = \omega_0\omega_1\omega_2\cdots &\mapsto \theta\omega := \omega_1\omega_2\omega_3\cdots
\end{align*}

Furthermore, let $E$ be a $d$-dimensional manifold and $M$ the compact closure of a $d$-dimensional submanifold of $\mathbb{R}^n$ and denote $\left\{\partial\right\} = E\backslash M$. Suppose
\begin{align*}
    f : K \times E &\rightarrow E\\
    (\omega_i,x) &\mapsto f_{\omega_i}(x)
\end{align*}
is continuously differentiable and $\mathrm{D} f$ is invertible on $K \times M$. Now define recursively the cocycle $(\theta, \varphi)$ as

$$\varphi_{n+1}(\omega, x)=
\begin{cases}
f_{\omega_n} \circ \varphi_n(\omega, x) & ,\varphi_n(\omega, x) \in M,\\
\partial & ,\varphi_n(\omega, x) = \partial.
\end{cases}
$$
Thus $(\theta, \varphi)$ forms a $\mathcal C^1$-random dynamical system with absorption at $\{\partial\}$ with linearised flow~$\Phi$ generated by $\Phi_1 = \mathrm{D} f$. Assume finally that the absorbed Markov process $(\varphi_n)_{n \in \mathbb{N}_0}$ fulfils \cite[ Hypothesis (H) \& Theorem 3-(M1)]{Castro2021ExistenceApproach} so that Hypothesis \ref{(H)} holds. Now by compactness of $K\times M$, $(\Theta,\Phi)$ fullfils the integrability condition \eqref{Qnu} and thus Theorems \ref{thm:fullFK} and \ref{thm:fullMET} and Corollary \ref{T2.12} apply and there exist conditioned Lyapunov exponents $\Lambda_1, \dots, \Lambda_d$ such that for $\rho^k$-almost every $(x, v) \in \mathrm{Gr}_k(M)$, $k \leq d$,
\begin{align*}
    \lim_{n\to\infty } \mathbb E_x \left[\left|\lambda^{(k)} - \frac{1}{n}\log \left\|\Bigwedge^k\Phi_n v\right\| \right| \bigg|~   \tau > n \right]  = 0.
\end{align*}
\end{example}

\subsection{Application to stochastic differential equations} \label{Sec2:SDE}

In this section, we aim to apply the results of Section \ref{CondLyapunovExp} to stochastic differential equations with escape. Let $M\subset \mathbb R^d,$ 
\begin{align*}
    (\Omega, (\mathcal{F}_t)_{t\geq 0}, \mathcal{F},\mathbb P) = \left(\mathcal C_0(\R_+,\R^m),\left(\sigma(\pi_s,0\leq s\leq t)\right)_{t\geq 0}, \sigma(\pi_s,s\geq 0),\mathbb P\right),
\end{align*}
where $\mathbb P$ is the Wiener measure, and consider the stochastic differential equation
\begin{align}
    \d X_t  = V_0 (X_t) \d t + \sum_{i=1}^{m} V_i(X_t) \circ \d W^i_t,\ X_0\in M, \label{SDE}
\end{align}on $M,$ where $(W_t^1,\ldots, W_t^m)$ denotes an $m$-dimensional standard Brownian motion, and ${V_i:M\to\mathbb R^d}$ are vector fields on $M$.

The assumption below ensures that the distribution of the process $(\varphi_t)_{t\geq 0}$, generated by the SDE \eqref{SDE}, converges exponentially to the quasi-stationary distribution, i.e.~it satisfies Hypothesis \ref{(H)}. %
\begingroup
\renewcommand\thehypothesis{(H\textsubscript{SDE})} 
\begin{hypothesis} \label{(B)} We say that the stochastic differential equation \eqref{SDE} fulfils Hypothesis $\mathrm{\ref{(B)}}$ if
\begin{enumerate}
    \item[${\mathrm{(H1_{\mathrm{SDE}})}}$] $M$ is an open connected and bounded subset of $\mathbb R^d$ with $\mathcal C^2$-boundary.
    \item[${\mathrm{(H2_{\mathrm{SDE}})}}$] The vector fields  $\{V_i:M\to\mathbb R^d\}_{i=0}^{m}$ admit  vector field extensions $\{\widetilde{V}_i:\mathbb R^d\to\mathbb R^d\}$ such that
    \begin{itemize}
        \item[(i)] $\left.\widetilde{V}_i\right|_{M} = V_i, \ \text{for every }i\in\{0,1,\ldots,m\};$ and

        \item[(ii)] $ \widetilde V_0\text{ is a }\mathcal C^1\text{-vector field}\ \text{and }\widetilde V_i\text{ is a {\revision smooth} vector field}$ for every  $i\in\{1,\ldots,m\}$. Note that since $\overline{M}$ is compact, this ensures the boundedness of the derivatives of $V_i$ on $M$.%
    \end{itemize}
    \item[${\mathrm{(H3_{\mathrm{SDE}})}}$] {\revision For every $x\in \partial M$, there exist an outward unit normal vector $v$ at $x$, and $i\in\{1,\ldots,m\}$ such that $\langle V_i(x), v\rangle \neq 0$.}

\item[${\mathrm{(H4_{\mathrm{SDE}})}}$] {\revision Equation \eqref{SDE} satisfies the strong Hörmander condition, i.e. for every $x\in M$, we have that $\mathrm{Span}([V_1,\ldots, V_m](x)) = T_x M \simeq  \mathbb R^d,$ where 
$$[V_1,\ldots V_m](x) := \left\{ v\in \mathbb R^d\ \middle\vert \begin{array}{l}
  \text{there exists }n\in\mathbb N\ \text{and }i_{1},i_2,\ldots, i_n\in\{1,\ldots,m\}\ \text{such that}\\  
\ v=[V_{i_1},[V_{i_2},[\ldots, [V_{i_{n-1}},V_{i_n}]\ldots ] ](x)\end{array}\right\}$$
and given two smooth vector fields $Y$ and $Z$ on $M$, we denote $[Y,Z]$ as the Lie bracket between $Y$ and $Z$.
}
\end{enumerate}
\end{hypothesis}
\endgroup

{\revision Observe that Hypothesis ${\mathrm{(H3_{\mathrm{SDE}})}}$ and ${\mathrm{(H4_{\mathrm{SDE}})}}$ are implied if the generator of (\ref{SDE}) is uniformly elliptic, i.e.~there exists $c>0$ such that for all $x\in M$ and for all $\xi \in \R^d\backslash\{0\}$
$$\sum_{i=1}^d\sum_{j=1}^d\sum_{k=1}^m \xi_i V_k^i(x) V_k^j(x) \xi_j > c \left\|\xi\right\|^2.$$
  }

    {\revision
Moreover, if uniform ellipticity is verified, the regularity of the vector fields $V_1,\ldots, V_m$ can be loosened to $\mathcal C^2$ rather than smooth.
}

The following theorem yields that the Lyapunov exponents given by the $Q$-process equate to the conditioned Lyapunov exponents conjectured in \cite{Engel2019ConditionedSystems}.

\begin{theorem}\label{T2.10}
Let $M$ be an open connected and bounded subset of $\mathbb R^d$ with $\mathcal C^2$-boundary, and suppose that the stochastic differential equation \begin{align}
    \d X_t  = V_0 (X_t) \d t + \sum_{i=1}^{m} V_i(X_t) \circ \d W^i_t,\ X_0\in M, \label{SDE2}
\end{align} 
satisfies Hypothesis $\mathrm{\ref{(B)}}.$ Then the random dynamical system $(\theta,\varphi)$ associated to (\ref{SDE2}) satisfies Hypothesis $\mathrm{\ref{(H)}}$. Furthermore, the linearised flow $\Phi_t$ fulfils condition \eqref{eqn:IC+} so that the multiplicative ergodic theorem holds.

Let $\nu$ denote the unique quasi-ergodic distribution of $(\theta,\varphi)$ on $M$ and $\sigma^k$ the Haar measure on the $k\textsuperscript{th}$ Grassmannian of $\R^d$,  $\mathrm{Gr}_k(T_x M) = \mathrm{Gr}_k(\R^d)$, defined in Remark \ref{remarkgrass}. Then for $(\nu\times\sigma^k)$-almost every $(x,v) \in \mathrm{Gr}_k(M) \simeq M\times \mathrm{Gr}_k(\R^d)$, we have that for every  $\varepsilon >0,$ 
\begin{align}
    \lim_{t\to\infty }\mathbb P_x \left[ \left|\lambda^{(k)} - \frac{1}{t}\log \left\|\Bigwedge^k\Phi_t v\right\| \right| >\varepsilon \bigg|~   \tau > t \right]  = 0,\label{C2}
\end{align}
where
$$\lambda^{(k)}=\Lambda_1+\cdots +\Lambda_k$$
and $\Lambda_1,\ldots, \Lambda_k$ are the $k$ first Lyapunov exponents given by Theorem \ref{thm:fullFK}.
\end{theorem}

Theorem \ref{T2.10} is proved in Section \ref{PT2.10}. This result generalises the theorem of convergence of finite-time Lyapunov exponents towards the average conditioned first Lyapunov exponent in conditional probability \cite[Theorem 3.9]{Engel2019ConditionedSystems}. Indeed, now recall that the Lyapunov exponents $\Lambda_i$ are no longer defined with respect to a quasi-ergodic distribution of an extended process. The existence and uniqueness of such a measure need not be assumed, and the $Q$-process setting covers cases where this does not hold (see Example \ref{ex:uncoupled}). However, in non-degenerate examples, one can in general compute the full spectrum of Lyapunov exponents with formulae we derive in Section \ref{sec:4}.

\begin{example}[Uncoupled stochastic differential equation]
    Consider the simple uncoupled two-dimensional SDE
    \begin{equation}
    \begin{cases}
        \d X_t &= (X_t - X_t^3)\d t + \sigma_1 \d W^1_t\\
        \d Y_t &= (Y_t - Y_t^3)\d t + \sigma_2 \d W^2_t
    \end{cases}
    \label{eqn:decoupled}
    \end{equation}
    with absorption at the boundary of the square domain $[-1.5, 1.5]\times [-1.5,1.5]$. Denote by $\varphi = (\varphi^1, \varphi^2)$ the generated random dynamical system and by $\nu_1$ and $\nu_2$ the unique quasi-ergodic distributions on $[-1.5,1.5]$ of $(\varphi^1_t)_{t\geq 0}$ and $(\varphi^2_t)_{t\geq 0}$ respectively. Hence, one can define
    
    $$\Lambda_i = \int_{-1.5}^{1.5} (1 - 3 z^2) \nu_i (\d z).$$
the conditioned average Lyapunov exponent achieved by $\varphi^i$. In this case, the conditioned process $(\varphi_t)_{t\geq 0}$ converges exponentially to quasi-stationarity and has quasi-ergodic distribution $\nu_1 \times \nu_2$. However, the top conditioned Lyapunov in the sense of Engel et al.~\cite{Engel2019ConditionedSystems} cannot be defined. Indeed the process $(\varphi_t, s_t)_{t\geq 0}$, where

$$s_t(\omega, x, v) :=\frac{D\varphi_t(\omega, x) v}{\left\|D\varphi_t(\omega, x) v\right\|} \qquad \mbox{for some }v\in \mathbb{S}^{1},$$
    does not have unique quasi-stationary and quasi-ergodic distributions. However, one can deduce that
    \begin{itemize}
        \item[Case 1:] $\sigma_1 > \sigma_2\implies \Lambda_1>\Lambda_2$: The conditioned process displays two (quasi)-ergodic components $\nu_1 \times \nu_2 \times \delta_{(\pm 1,0)}$ and $\nu_1 \times \nu_2 \times \delta_{(0,\pm1)}$ achieving $\Lambda_1$ and $\Lambda_2$ respectively. Thus, this system yields a Lyapunov spectrum $\{(\Lambda_1, 1),(\Lambda_2, 1)\}$. In addition, with respect to the Lebesgue measure, almost every $v \in \mathbb{S}^{1}$ achieves $\Lambda_1$. One can be even more precise: all $v \in \mathbb{S}^1\backslash\{(0, \pm 1)\}$ achieve $\Lambda_1$, and $(0, \pm 1)$ achieves $\Lambda_2$. In other words, this exhibits the structure of an Oseledets flag $\{(0,0)\}\subset \mathrm{Span}\{(0,1)\} \subset \R^2$.
        \item[Case 2:] $\sigma_1 = \sigma_2\implies \Lambda_1 = \Lambda_2$: In this case, we obtain that every $v \in \mathbb{S}^{1}$ achieves $\Lambda_1 = \Lambda_2$. Thus the Lyapunov spectrum is $\{(\Lambda_1, 2)\}$, i.e.~$\Lambda_1$ has multiplicity $2$. Here, the Oseledets flag is given by $\{0\} \subset \R^2$.
    \end{itemize}
    We emphasise that Theorem \ref{T2.10} covers such degenerate examples.
    \label{ex:uncoupled}

\end{example}

%
%
%
%
%
%

%

%
%

%

%
%
%
%
   
%
%
%
%
%
%
%
%
    
%
%
%
    
%
%
    
%
    
%
%
%

%
%
%
%
%
%
%

%

%
%
%
%

%

%
%
%
%
%
%
%
%
%
%

%
%
%
%
%
%
%
%
%
%

%
%
%
%
%
%
%
%
%

%
%
%
%
%
%
%
%
%
%

%
%
%
%
%
%
%
 
%
%

%

%
%
%
\section{The Q-process for random dynamical systems}
\label{sec:3}

In this section, we prove the existence of the $Q$-process process for a random dynamical system fulfilling Hypothesis \ref{(H)} and the ergodicity of $\mathbb Q_\nu$ under $\Theta$. In this paper, we  restrict  to the two following noise spaces:
 
\begin{enumerate}
     \item[$(i)$]
In the case $\mathbb T = \mathbb R_+$, we always consider
\begin{align}
    (\Omega, (\mathcal{F}_t)_{t\geq 0}, \mathcal{F}) = \left(\Omega,\left(\sigma(\pi_s,0\leq s\leq t)\right)_{t\geq 0} , \sigma(\pi_s,s\geq 0)\right),\label{l1}
\end{align}
where $\Omega \in \left\{  \mathcal{D}(\R_+,\R^m), \mathcal{C}_0(\R_+,\R^m)\right\}$ and
\begin{align*}
    \mathcal C_0(\R_{+},\R^m) &=\{\omega:\R_+\to \R^m; \ \omega\text{ is a  continuous and }\omega(0) = 0 \},\\
    \mathcal D(\R_{+},\R^m) &=\{\omega:\R_+\to \R^m; \; \ \omega\text{ is \textit{càdlàg}}\},
\end{align*}

and
$$\pi_t : \omega\in \Omega \to \omega(t)\in \R^m. $$
\item[$(ii)$] In the case that $\mathbb T = \mathbb N_0$,
\begin{align}
    (\Omega, (\mathcal{F}_n)_{n\in \mathbb{N}_0}, \mathcal{F}) = \left( X^{\mathbb N_0},\left(\sigma(\pi_m,0\leq m\leq n)\right)_{n\in\mathbb N_0}, \sigma(\pi_m,m\geq 0)\right),\label{l2}
\end{align}
where $X$ is a Polish space,
$$ X^{\mathbb N_0} = \{f:\mathbb N_0 \to X\}, $$
and 
$$\pi_n: \omega \in X^{\mathbb N_0} \to \omega(n) \in X. $$
\end{enumerate}

Let us also recall some properties of the conditioned process under quasi-stationarity.

\begin{proposition}
\label{prop:QSD_QED}
If Hypothesis \ref{(H)} is fulfilled, then \label{Prop25}
\begin{enumerate}
    \item[(i)] $\nu(\d x) = \eta(x) \mu(\d x);$
    \item[(ii)] $\int_M \mathcal P^t(x,\cdot) \d \mu = e^{-\beta  t}\mu(\cdot)$ for every $t\in\mathbb T$; 
    \item[(iii)] $\int_M \mathcal P^t(\eta)(x)\mu(\d x) = e^{-\beta  t}$ for every $t\in\mathbb T.$\label{prop:A}
\end{enumerate}
\end{proposition}

Proposition \ref{prop:A} is proved in Appendix \ref{Appendix:B}.

\subsection{{\revision Proof of Proposition \ref{Q-process}}}
\label{subsec:Q-process}
In the following, we prove Proposition \ref{Q-process}. The proof relies on the following proposition.

\begin{proposition}
Let $(\theta,\varphi)$ be an absorbed random dynamical system fulfilling Hypothesis $ \mathrm{\ref{(H)}} .$ Then, for every $x\in M,$ there exists a unique probability measure $\mathbb Q_x$ on  $(\Omega\times M,\mathcal G = \mathcal F \otimes \mathcal B(M)),$ such that for every $s\geq 0$ and $A\in\mathcal G_s,$
$$\mathbb Q_x[A] = \frac{e^{\beta  s}}{\eta(x)}\mathbb E_x[\mathbbm 1_A \eta\circ\varphi_s ] = \lim_{t\to\infty}\mathbb  P_x[A \mid \tau >t ].$$\label{marathon}
\end{proposition}

\begin{proof}

Let us fix $x\in M.$ We divide the proof in four steps.
\begin{step}{1} \it We show that for every $x\in M,$ $s\in \mathbb T$ and $A\in\mathcal G_s,$
$$\lim_{t\to\infty} \mathbb P[A\mid \tau(\cdot,x) > t] = \frac{e^{\beta  s}}{\eta(x)}\mathbb E\left[\mathbbm 1_A \ \eta(\varphi_s(\cdot,x))\right].$$
\end{step}
By a direct computation,

\begin{align*}
    \lim_{t \to\infty}\mathbb P_x [A\mid \tau >t] &\stackrel{\phantom{\mathrm{(H3)}}}{=} \lim_{t \to\infty} \frac{\mathbb E_x [\mathbbm 1_A \mathbbm 1_{\{\tau >t\}}]}{\mathcal P^t(x,M)}\\
    &\stackrel{\phantom{\mathrm{(H3)}}}{=}\lim_{t \to\infty} \frac{\mathbb E_x\left[\mathbb E_x [\mathbbm 1_A \mathbbm 1_{\{\tau >t\}}\mid \mathcal G_s] \right]}{\mathcal P^t(x,M)}\\
    &\stackrel{\phantom{\mathrm{(H3)}}}{=}\lim_{t \to\infty}\frac{\mathbb E_x\left[\mathbbm 1_A \mathcal P^{t-s}\left(\varphi_{s},M \right)\right]}{\mathcal P^t(x,M)}\\
    &\stackrel{\phantom{\mathrm{(H3)}}}{=}\lim_{t\to\infty}\frac{e^{\beta  s}\mathbb E_x\left[\mathbbm 1_A e^{\beta  (t-s)}\mathcal P^{t-s}\left(\varphi_{s},M \right)\right]}{e^{\beta  t}\mathcal P^t(x,M)}\\
    & \stackrel{\mathrm{(H3)}}{=}\frac{e^{\beta  s}}{\eta(x)} \mathbb E[\mathbbm 1_A \eta(\varphi_s(\cdot,x))] .
\end{align*}
This proves Step 1.

\begin{step}{2} We show that there exists a measure $\widetilde{\mathbb{Q}}_x$ on $(\Omega ,\mathcal F ),$ such that
$$\left.\widetilde{\mathbb{Q}}_x\right|_{{\mathcal F}_t} [\d \omega] = \frac{e^{\beta  t}}{\eta(x)}\eta (\varphi_t(\omega,x)) \mathbb P[\d \omega]. $$

\end{step}

In the case $\mathbb T = \mathbb N_0$ this follows immediately from Kolmogorov's extension theorem (see \cite[Theorem II 3.26.1]{Rogers2000DiffusionsMartingales}). In the case $\mathbb T = \mathbb R_+$ and $\Omega =\mathcal D(\R_+,\R^m)$ consider the measure space isomorphism 
\begin{align*}
    \mathcal{E}:  \left(\mathcal D(\R_+,\R^m), \sigma\left( (\pi_s)_{s\geq 0}\right) \right) &\rightarrow \left(\widehat{\Omega}, \widehat{\mathcal{F}}\right) \\
    \omega &\mapsto \big(\omega|_{[0,1]}, \omega|_{[1,2]}-\omega(1), \omega|_{[2,3]}- \omega(2),\dots\big),
\end{align*}
where  $$(\widehat{\Omega}, \widehat{\mathcal{F}}) :=\left(\mathcal D([0,1],\R^m), \sigma\left(\left(\pi_s\right)_{0 \leq s\leq 1}\right)\right)\otimes (\mathcal D_0([0,1],\R^m), \sigma((\pi_s)_{0 \leq s\leq 1})^{\otimes \mathbb N_0}$$ and $\mathcal D_0([0,1],\R^m) =\{\omega \in \mathcal D([0,1],\R ^m); \omega(0) = 0\}.$

Consider the filtration $(\widehat{\mathcal F}_n)_{n\in\mathbb N}$ on $(\widehat{\Omega}, \widehat{\mathcal{F}})$ by $\widehat{\mathcal F}_n := \mathcal E (\mathcal F_n)$. Note that for each $n \in \mathbb N$ the $\sigma$-algebra $\widehat{\mathcal F}_n$ is precisely the one induced by projection on the first $n$ components of $\mathcal D([0,1],\R^m)\otimes \mathcal D_0([0,1],\R^m)^{\otimes \mathbb N_0}$. Since $\mathcal D([0,1],\R^m)$ and  $\mathcal D_0([0,1],\R^m)$ are a Polish spaces, we can apply Kolmogorov's extension theorem to the measures $$\widehat{\mathbb Q}_x^n := \mathcal{E}_*\left(\frac{e^{\beta  t}}{\eta(x)}\eta (\varphi_t(\omega,x) \mathbb P[\d \omega]\right)\text{ on }\left( \mathcal D([0,1],\R^m)\times \mathcal D_0([0,1],\R^m)^{\times  \mathbb N_0},\widehat{\mathcal F}_n\right) \ \text{for every }n\in\mathbb N_0,$$ to get a measure $\widehat{\mathbb Q}_x$ on $\big(\widehat{\Omega}, \widehat{\mathcal{F}}\big)$ with $\widehat{\mathbb Q}_x^n = \widehat{\mathbb Q}_x|_{\widehat{\mathcal F}_n}$. Now we can set $\widetilde{ \mathbb Q}_x := \mathcal{E}^{-1}_* \widehat{\mathbb Q}_x$ to get the desired measure. 

If $\Omega = \mathcal C_0(\mathbb R_+, \mathbb R^m)$ the exact same argument can be applied changing $\mathcal D$ to  $\mathcal C_0$ in the above proof. 

\begin{step}{3} We complete the proof.
\end{step}

Consider $\mathbb Q_x := \widetilde{\mathbb Q}_x\times \delta_x.$ From Steps 1--2 it is clear that $\mathbb Q_x$ is a Borel measure on ${\mathcal G = \mathcal F\otimes \mathcal B(M)}$ and for every $x\in M,$  $s\in \mathbb T$ and $A\in\mathcal G_s,$
$$\mathbb Q_x[A] = \frac{e^{\beta  s}}{\eta(x)}\mathbb E_x[\mathbbm 1_A \eta\circ\varphi_s ] = \lim_{t\to\infty}\mathbb  P_x[A \mid \tau >t ].$$
The uniqueness of $\mathbb Q_x$ follows directly from the monotone class theorem and
$$\sigma\left(\bigcup_{s\geq 0} \mathcal G_s\right) = \mathcal G. $$
This finishes the proof.
\end{proof}

In the following, we prove a useful lemma.
\begin{lemma}\label{Qmarkovlemma}
Let $A_1 \in \mathcal G_t$ and $A_2 \in \mathcal G_s$. Then we have for all $x \in \mathcal M$
$$\mathbb Q_x\big(A_1 \cap \Theta_t^{-1}(A_2)\big) = \int_{\mathcal M} \mathbb Q_y (A_2)~\mathbb Q_x\big(A_1 \cap \{\varphi_t \in dy\}\big).$$
As an easy consequence, we also have
$$\mathbb Q_\nu\big(A_1 \cap \Theta_t^{-1}(A_2)\big) = \int_{\mathcal M} \mathbb Q_y (A_2)~\mathbb Q_\nu\big(A_1 \cap \{\varphi_t \in dy\} \big).$$
\end{lemma}

\begin{proof}
Clearly, the second equality follows from the first by integration with respect to $\nu$. Thus, we only show the first equality. Observe that
\begin{align}
   \mathbb Q_x\big(A_1 \cap \Theta_t^{-1}(A_2)\big) &=\mathbb{E}_{x}^{\mathbb Q}\left[\mathbbm 1_{A_1} \cdot \mathbbm 1_{A_2} \circ \Theta_t\right]\nonumber\\
   &= \frac{e^{\beta (t+s)}}{\eta(x)} \mathbb E_x\left[\mathbbm 1_{A_1} \cdot \mathbbm 1_{A_2}\circ \Theta_t \cdot \eta\circ\varphi_{s+t} \right]\nonumber\\
    &=\frac{e^{\beta (t+s)}}{\eta(x)} \mathbb E_x\left[\mathbbm 1_{A_1} \cdot \mathbbm 1_{A_2}\circ \Theta_t \cdot\eta\circ \varphi_{s} \circ \Theta_t  \right].\label{D1}
\end{align}
Since the $\sigma$-algebras $\mathcal F_t$ and $\theta_{t}^{-1}\mathcal F_s$ are $\mathbb P$-independent, we obtain
\begin{align}
 \mathbb E_x\left[\mathbbm 1_{A_1} \cdot \mathbbm 1_{A_2}\circ \Theta_t \cdot\eta\circ \varphi_{s} \circ \Theta_t  \right]&=  \mathbb E_x\left[\mathbbm 1_{A_1}\cdot \mathbb E\left[\mathbbm 1_{A_2}\circ (\theta_t,\varphi_t(\cdot,x)) \cdot\eta\circ \varphi_{s} \circ (\theta_t,\varphi_t(\cdot,x))  \mid \mathcal  F_t\right]\right]\nonumber\\
 &= \mathbb E_x\left[\mathbbm 1_{A_1}\cdot \mathbb E_{\varphi_t(\cdot,x)}\left[\mathbbm 1_{A_2}\circ (\theta_t,\cdot) \cdot\eta\circ \varphi_{s} \circ (\theta_t,\cdot)  \mid \mathcal  G_t\right]\right]\nonumber\\
  &= \mathbb E_x\left[\mathbbm 1_{A_1}\cdot \mathbb E_{\varphi_t(\cdot,x)}\left[\mathbbm 1_{A_2} \cdot\eta\circ \varphi_{s}  \right]\right]\label{D2}.
 \end{align}
Combining (\ref{D1}) and (\ref{D2}), we achieve
\begin{align*}
    \mathbb Q_x\big(A_1 \cap \Theta_t^{-1}(A_2)\big) &= \frac{e^{\beta t}}{\eta(x)} \mathbb E_x\left[\mathbbm 1_{A_1}\cdot \eta \circ\varphi_t \cdot \frac{e^{\beta s}}{\eta \circ \varphi_t}\mathbb E_{\varphi_t(\cdot, \cdot)}\left[\mathbbm 1_{A_2} \cdot\eta\circ \varphi_{s}\right]\right]\\
    &= \frac{e^{\beta t}}{\eta(x)} \mathbb E_x\left[\mathbbm 1_{A_1}\cdot   \eta \circ\varphi_t \cdot \mathbb Q_{\varphi_t(\cdot, \cdot)}[A_2]\right]\\
    &= \mathbb{E}_{x}^{\mathbb Q}\left[\mathbbm 1_{A_1} \cdot \mathbb Q_{\varphi_t(\cdot, \cdot)}[A_2]\right] \\
    &= \int_{\mathcal M} \mathbb Q_y (A_2)~\mathbb Q_x\big(A_1 \cap \{\varphi_t \in dy\}\big),
\end{align*}
which yields the statement.
\end{proof}

Now we can prove Proposition \ref{Q-process}.
\begin{proof}[Proof of Proposition \ref{Q-process}]
Let $\{\mathbb Q_x\}_{x\in M}$ be the family of measures given by Proposition \ref{marathon}. We divide the proof into three steps.

\begin{step}{1} \it We show that $\{\mathcal Q^t\}_{t\in\mathbb T}$ fulfils the Chapman--Kolmogorov equation.
\end{step}
Given $t,s\in\mathbb T$ and $B\in\mathcal B(M),$ we have that 
\begin{align*}
    \mathcal Q^{t+s}(x,B) &= \frac{e^{ \beta  (t+s) }}{\eta(x)} \mathcal P^{t+s}(\mathbbm 1_B \eta)(x)= \frac{e^{\beta  t}}{\eta(x)} \mathcal P^{t}\left(\frac{e^{\beta  s}}{\eta(\cdot)}\mathcal P^s(\mathbbm 1_B \eta )(\cdot) \eta(\cdot)\right)(x)\\
    &= \frac{e^{\beta  t}}{\eta(x)} \mathcal P^t\left( \mathcal Q^s(\mathbbm 1_B) (\cdot) \eta(\cdot) \right)(x)\\
    &= \int_{M} \mathcal Q^s(y,B) \mathcal Q^t(x,\d y). 
\end{align*}
\begin{step}{2}
We show that
\begin{align}
    \left(\Omega\times M, \left(\mathcal{G}_t\right)_{t\in\mathbb T}, \left(\varphi_t\right)_{t\in\mathbb T}, \left(\mathcal Q^t\right)_{t\in\mathbb T}, (\mathbb Q_x)_{x \in M}\right)\label{Q}
\end{align}
is a Markov Process. 
\end{step}

To conclude that (\ref{Q}) is a Markov process, the only non-trivial property that must be verified is that for all $t,s\in \mathbb T$ and  every bounded measurable function $f$ on $M$
$$\mathbb{E}_{x}^{\mathbb Q}[f\circ \varphi_{t+s} \mid \mathcal G_s] = \mathcal (Q^{t}f)(\varphi_s)\quad \mbox{$\mathbb Q_x$-almost surely.} $$
To verify this let $A\in\mathcal G_s.$ We can compute
\begin{align*}
    \mathbb{E}_{x}^{\mathbb Q}[\mathbbm 1_A \cdot f\circ \varphi_{t+s}] \stackrel{\phantom{\mathrm{Lem.\  \ref{Qmarkovlemma}} }}{=} & \mathbb{E}_{x}^{\mathbb Q}[\mathbbm 1_A \cdot f\circ \varphi_{t}\circ \Theta_s]\\
   \stackrel{\mathrm{Lem.\  \ref{Qmarkovlemma}} }{=}& \int_M \mathbb E_{\mathbb Q_y} [f\circ \varphi_{t}]~\mathbb Q_x\big(A \cap \{\varphi_s \in dy\}\big) \\
    \stackrel{\phantom{\mathrm{Lem.\  \ref{Qmarkovlemma}} }}{=} & \int_M \mathcal (Q_t f)(y)~\mathbb Q_x\big(A \cap \{\varphi_s \in dy\}\big) \\
    \stackrel{\phantom{\mathrm{Lem.\  \ref{Qmarkovlemma}} }}{=} & \mathbb{E}_{x}^{\mathbb Q} [\mathbbm 1_A \cdot \mathcal (Q_t f)(\varphi_s)].
\end{align*}

\begin{step}{3}\it  We show that $(\ref{Q})$ admits the quasi-ergodic distribution $\nu$ as its unique ergodic invariant measure and for every $x\in M$
$$\lim_{t\to\infty}\|\mathbb Q_x (\varphi_t\in \cdot) - \nu(\cdot)\|_{TV}=0. $$

\end{step}

Let $B\in\mathcal B(M),$ then we have that for every $t\in\mathbb T,$
\begin{align}
    \left|\mathbb Q_x [\varphi_t \in B] - \nu(B)\right| =& \left| \frac{e^{\beta  t}}{\eta(x)}\mathbb E_x [\mathbbm 1_B \circ \varphi_t \cdot \eta\circ \varphi_t ] - \nu(B)\right| \nonumber\\
    =&  \left| \frac{e^{\beta  t} \mathbb P_x[\tau>t]}{\eta(x)}\mathbb E_x [\mathbbm 1_B \circ \varphi_t \cdot \eta\circ \varphi_t\mid \tau >t ] - \int_B \eta \d \mu\right| \nonumber\\
    \leq& \left| \frac{e^{\beta  t} \mathbb P_x[\tau>t]}{\eta(x)} -1\right|\|\eta\|_{\infty}  +\left|  \mathbb E_x [ \mathbbm 1_B\circ \varphi_t\cdot \eta\circ \varphi_t\mid \tau >t ] - \int_B \eta \d \mu\right|\label{rei}. 
\end{align}

From $ \mathrm{(H2)}$--$\mathrm{(H3)}$ we have that 
\begin{align}
    \left|  \mathbb P_x [ \mathbbm 1_B\circ \varphi_t\cdot \eta\circ \varphi_t\mid \tau >t ] - \int_B \eta \d \mu\right| \leq \|\eta\|_{\infty} C(x) e^{-\alpha t} \label{rei1}
\end{align}
and
\begin{align}
    \lim_{t\to\infty }\left| \frac{e^{\beta  t} \mathbb P_x[\tau>t]}{\eta(x)} -1\right| = 0.\label{rei2}
\end{align}

From equations $(\ref{rei}),$ $(\ref{rei1})$ and $(\ref{rei2}),$ we obtain
$$\lim_{t\to\infty}\left\|\mathbb Q_x\left[\varphi_t \in \cdot~\right] - \nu\right\|_{TV} =0.$$

The above equation implies that $\nu$ is the unique stationary and therefore ergodic measure of the Markov process $\left(\Omega\times M, \left(\mathcal{G}_t\right)_{t\in\mathbb T}, \left(\varphi_t\right)_{t\in\mathbb T}, \left(\mathcal Q^t\right)_{t\in\mathbb T}, \left(\mathbb Q_x\right)_{x \in M}\right)$.
\end{proof}

\begin{remark}
Note that the initial probabilities $(\mathbb Q_x)_{x\in M}$ associated to the $Q$-process of the stochastic differential equation $(\ref{SDE2})$ do not depend on chosen modifications of $\varphi$ (see \cite[Definition~II. 36.2]{Rogers2000DiffusionsMartingales}). Indeed, let $\widetilde \varphi$ be a modification of $\varphi$ and consider the stopping time
$$\widetilde \tau(\omega,x) = \inf \{t \geq 0:\ \widetilde \varphi (t,\omega,x) \not \in M\}.$$
We have that for every $x\in M$
$$\tau  = \widetilde \tau\ \qquad \mbox{$\mathbb P_x$-almost surely}. $$
This implies that for every $t\geq 0,$
$$\mathbb P_x [\ \cdot \mid \widetilde \tau >t ] = \mathbb P_x [\ \cdot \mid \tau >t ],$$
and therefore $\widetilde \varphi$ and $\varphi$  generate the same family of probabilities $(\mathbb Q_x)_{x\in M}$.  Furthermore, properties such as continuity, \textit{càdlàg} paths and the cocycle property of $\varphi$ under $\mathbb{P}_x$ are preserved under $\mathbb{Q}_x$.
\label{rmk:modification}
\end{remark}

\subsection{The Q-process dynamical framework {\revision and proof of Theorem~\ref{thm:IntroErgMeas}}}
\label{subsec:proof_TheoremA}
Let $(\theta,\varphi)$ be a random dynamical system fulfilling Hypothesis $ \mathrm{\ref{(H)}} $, $(\mathbb Q_x)_{x\in M}$  the family of measures given by Proposition \ref{marathon} and $\nu$ the unique quasi-ergodic distribution of $(\theta,\varphi)$ on $M$ given by Hypothesis $\mathrm{(H1)}$.

In this section we prove that the measure $$\mathbb Q_\nu := \int_M \mathbb Q_x \nu(\d x),$$
on $\Omega\times M$ is an ergodic measure for the skew product $(\Theta_t)_{t\in\mathbb T}$ of the random dynamical system~$(\theta,\varphi)$. 

\begin{lemma}
Let $(\theta,\varphi)$ be a random dynamical system fulfilling Hypothesis~$\mathrm{\ref{(H)}}$, and let $(\Theta_t)_{t\in\mathbb T}$ be its skew product. Then 
$$\left(\mathbb P \times \mu\right)\left( \Theta_t^{-1}(C)\right) =  e^{-\beta  t}\left(\mathbb P \times \mu\right) (C), \ \forall\ t\geq 0 \ \text{and }C\in \mathcal G,$$
where $$\Theta_t(\omega,x) = \left(\theta_t (\omega), \varphi_t(\omega,x) \right).$$\label{lemito}
\end{lemma}
\begin{proof}
Fix $s\geq 0$. Consider $A\times B \in \mathcal F_s \otimes \mathcal B(M),$ then
\begin{align*}
   \left(\mathbb P\times \mu \right)\left(\Theta_t^{-1}(A\times B)\right)&=\int_{\Omega\times M} \mathbbm 1_{A\times B} (\Theta_t (\omega,x)) \left(\mathbb P\times \mu  \right)(\d \omega\times \d x)\\
   &=  \int_{\Omega\times M} \mathbbm 1_{A\times B} (\theta_t (\omega),\varphi_t(\omega,x)) \left(\mathbb P\times \mu  \right)(\d \omega\times \d x)\\\
    &=\int_{ M} \mathbb E\left[\mathbbm 1_A\circ \theta_t(\cdot)\ \mathbbm 1_B \circ \varphi_t(\cdot,x)  \right]\mu(\d x)\\
    &= \int_M \mathbb P[A] \mathbb E[\mathbbm 1_B \circ \varphi_t(\cdot,x)]\mu(\d x)\\
    &= \mathbb P[A] \int_M  \mathcal P^t(x,B)\mu(\d x) = e^{-\beta  t} \left(\mathbb P\times \mu \right)\left(A\times B\right).
\end{align*}

Since  $$\sigma\left( \left\{A\times B;\ A\in \bigcup_{s\geq 0}\mathcal F_s,\ B\in\mathcal B(M)\right\}\right) =\mathcal G,$$ we obtain from the monotone class theorem that 
$$\left(\mathbb P \times \mu\right) ( \Theta_t^{-1}(C)) =  e^{-\beta  t}\left(\mathbb P \times \mu\right) (C), \ \forall\ t\geq 0 \ \text{and }C\in \mathcal G,$$
which finishes the proof.
\end{proof}

In the following, we prove {\revision Theorem \ref{thm:IntroErgMeas}}.

\begin{proof}[Proof of {\revision Theorem \ref{thm:IntroErgMeas}}]
We show that the measure $\mathbb Q_\nu$ is strongly mixing under $\Theta$, i.e.~we have
\begin{equation}\label{mixingEquation}
    \lim_{t \rightarrow \infty} \mathbb{Q}_\nu \big(A_1 \cap \Theta_t^{-1}(A_2)\big) = \mathbb Q_\nu (A_1) \mathbb Q_\nu (A_2), 
\end{equation}
for all $A_1, A_2 \in \mathcal{G}$. In particular, this implies both ergodicity and invariance.

{\revision From {\cite[Theorem 1.17]{Walters1982AnTheory}} it is sufficient to 
 prove that (\ref{mixingEquation}) holds in the case that $A_1 \in \mathcal G_s$ and $A_2 \in \mathcal G_r$.}

Under these additional assumptions, we have
\begin{align*}
    \lim_{t\rightarrow \infty}\mathbb Q_\nu\big(A_1 \cap \Theta_t^{-1}(A_2)\big) &\stackrel{}{\mathmakebox[\widthof{$\stackrel{\mathrm{Dom.\ Conv.\ Thm. }}{=}$}]{=}} \lim_{t\rightarrow \infty}\mathbb Q_\nu\big(A_1 \cap \Theta_{s}^{-1}\Theta^{-1}_{t-s}(A_2)\big)\\
    &\stackrel{ \mathrm{Lem.\  \ref{Qmarkovlemma} }}{\mathmakebox[\widthof{$\stackrel{\mathrm{Dom.\ Conv.\ Thm. }}{=}$}]{=}} \lim_{t\rightarrow \infty} \int_M \mathbb Q_x\big(\Theta^{-1}_{(t-s)}(A_2)\big)~\mathbb Q_\nu\big(A_1\cap \{\varphi_s \in dx\}\big)\\
   &\stackrel{ \mathrm{Dom.\ Conv.\ Thm. }}{=}   \int_M \left(\lim_{t\rightarrow \infty}\mathbb Q_x\big(\Theta^{-1}_{(t-s)}(A_2)\big)\right)~\mathbb Q_\nu\big(A_1\cap \{\varphi_s \in dx\}\big)\\
      &\stackrel{ \mathrm{Lem.\  \ref{Qmarkovlemma} }}{\mathmakebox[\widthof{$\stackrel{\mathrm{Dom.\ Conv.\ Thm. }}{=}$}]{=}} \int_M \left(\lim_{t\rightarrow \infty}\int_M \mathbb Q_y(A_2)~\mathbb Q_x\big(\varphi_{t-s}(\cdot, x) \in dy\big)\right)~\mathbb Q_\nu\big(A_1\cap \{\varphi_s \in dx\}\big)\\
      &\stackrel{ \mathrm{Prop.\ \ref{Q-process}}}{\mathmakebox[\widthof{$\stackrel{\mathrm{Dom.\ Conv.\ Thm. }}{=}$}]{=}}  \int_M \left(\int_M \mathbb Q_y(A_2)~\nu(dy)\right)~\mathbb Q_\nu\big(A_1\cap \{\varphi_s \in dx\}\big)\\
     &\stackrel{}{\mathmakebox[\widthof{$\stackrel{\mathrm{Dom.\ Conv.\ Thm. }}{=}$}]{=}}\mathbb Q_\nu(A_2)\int_M  ~\mathbb Q_\nu\big(A_1\cap \{\varphi_s \in dx\}\big)\\
    &\stackrel{}{\mathmakebox[\widthof{$\stackrel{\mathrm{Dom.\ Conv.\ Thm. }}{=}$}]{=}}\mathbb Q_\nu(A_1)\mathbb Q_\nu(A_2),
\end{align*} 
{\revision which concludes the proof of the Proposition.}

%
\end{proof}

\subsection{The multiplicative ergodic theorem for the Q-process}
\label{QMET}
\label{PT2.12}
\begin{proof}[Proof of Theorems~\ref{thm:fullFK} and \ref{thm:fullMET}]\label{METQP}
    We adapt the proof of \cite[Theorem 4.2.6]{Arnold1998RandomSystems} to one-sided time. 
    For more details, we refer to \cite[Chapters 3 \& 4]{Arnold1998RandomSystems}. Let $\langle \cdot, \cdot \rangle_{x}$ denote the Riemannian structure of $M$. Then there exists a global trivialisation $\psi : TM \rightarrow M\times \R^d$ (here in the sense of a bimeasurable bijection) such that for all $x \in M$
    $$\psi_x = \psi_{|T_xM} : (T_xM, \langle\cdot,\cdot\rangle_x)\rightarrow (\R^d, \langle\cdot, \cdot \rangle)$$
    is an isometry (see \cite[Lemma 4.2.4]{Arnold1998RandomSystems}). Furthermore, let
    \begin{align*}
        \widetilde{\Phi}_{t}((\omega, x), v) : \mathbb{T}\times (\Omega \times M) \times \R^d &\to \R^d\\
        (t, (\omega, x), v) &\mapsto\widetilde{\Phi}_{t}((\omega, x), v) := \psi_{\varphi_{t}(\omega, x)}\circ \Phi_t(\omega, x)\circ \psi_x^{-1}(v).
    \end{align*}
    Then $\widetilde{\Phi}$ forms a linear cocycle over the ergodic DS $\left(\Omega\times X, \mathcal{F}\otimes \mathcal{B}(M), \mathbb{Q}_{\nu}, (\Theta_t)_{t\in T}\right)$. Hence, we can apply the one-sided time versions of Furstenberg--Kesten theorem~\cite[Theorem 3.3.3]{Arnold1998RandomSystems} and the multiplicative ergodic theorem~\cite[Theorem 3.4.1]{Arnold1998RandomSystems} for linear cocycles on $\widetilde{\Phi}$. Since $\widetilde{\Phi}$ is Lyapunov cohomologous to $\Phi$, this yields the desired result.
    
    Observe that we can choose $\Delta$ as a subset of $\Xi$, since $\Xi$ is a full measure $\Theta_t$-forward-invariant set.
\end{proof}

Another characterisation of the Lyapunov exponents is obtained via the growth rates of $k$-volume forms along the trajectories of $\varphi$. In more detail, for $k\leq d$, $(\omega, x)\in \Delta$, and ${v_1, \ldots, v_k\in T_x M}$ linearly independent,
$$\lambda^{(k)}:=\Lambda_1 + \cdots +\Lambda_k = \lim_{t\to\infty}\frac{1}{t}\log\mathrm{Vol}\left(\Phi_t(\omega, x)v_1, \dots, \Phi_t(\omega, x)v_k\right)$$
where $\mathrm{Vol}(u_1, \ldots, u_k)$ denotes the volume of the parallelepiped spanned by $u_1, \ldots, u_k$.

\noindent For ${v = v_1 \wedge \cdots \wedge v_k \in \Bigwedge^k_0(T_x M)}$, one can compute this volume as
$$\left\|\Bigwedge^k\Phi_t(\omega,x) v\right\|:=\left\|\Phi_t(\omega,x) v_1 \wedge \cdots \wedge \Phi_t(\omega,x) v_k \right\| = \mathrm{Vol}\left(\Phi_t(\omega, x)v_1, \dots, \Phi_t(\omega, x)v_k\right)$$
where $\|\cdot \|$ is induced by the inner product on $\Bigwedge^k(T_x M)$ defined on $\Bigwedge^k_0(T_x M)$ as
$$\langle u_1\wedge\cdots\wedge u_k,v_1\wedge\cdots\wedge v_k\rangle:=\det (\langle u_i,v_j\rangle_x)_{k\times k.}$$

This motivates the introduction of
\begin{align}
    r^k_t(\omega, x, v) := \left\|\Bigwedge^k\Phi_t(\omega, x) v\right\|\qquad \text{and} \qquad s_t^k(\omega, x, v) = \frac{\Bigwedge^k\Phi_t(\omega, x) v}{\left\|\Bigwedge^k\Phi_t(\omega, x) v\right\|}\in \mathrm{Gr}_k(T_x M),\label{rands}
\end{align}
where $\mathrm{Gr}_k(T_x M)$ denotes the Grassmannian manifold $\mathbf{P}(\Bigwedge^k_0 T_x M)$, which is a $k(d-k)$-dimensional submanifold of the projective space $\mathbf{P}(\Bigwedge^k T_x M)$ (Plücker embedding \cite[Page 209]{Griffiths1994PrinciplesGeometry}). $\mathrm{Gr}_k(T_x M)$ can also be thought as the space of $k$-dimensional subspaces of $T_x M$. Note that, here we have implicitely identified antipodal points as they achieve the same Lyapunov exponents. {\revision We equip $\mathrm{Gr}_k(T_xM)$ with its Haar measure $\sigma^k_x$, i.e. which is invariant under rotations of $T_x(M)$ (See \cite[p. 325]{Baxendale1986TheDiffeomorphisms} for a construction).}

Denote the $k^{\mathrm{th}}$ Grassmannian bundle $\mathrm{Gr}_k(M)$ as the fiber bundle whose fiber at $x\in M$ is $\mathrm{Gr}_k(T_x M)$, i.e.
 $$\mathrm{Gr}_k(M) := \bigcup_{x \in M} \mathrm{Gr}_k(T_x M).$$

\begin{remark}
   In the case that $M$ is an open subset of $\mathbb R^m,$ we have that $\mathrm{Gr}_k(M) \cong M\times  \mathrm{Gr}_k(\R^d)$ as fiber bundle, which implies that every $\{\sigma_x^k\}_{x\in M}$ can be canonically identified with a single measure~$\sigma^k.$   \label{remarkgrass}
\end{remark}

This gives us the following corollary as a consequence of the multiplicative ergodic theorem.
\begin{corollary}\label{BaxendaleResult}
    For $k\leq d$, let $\Lambda_1, \ldots, \Lambda_k$ be as above and let $\sigma_x^k$ be the Haar measure on $\mathrm{Gr}_k(T_x M)$ above. Then for $\mathbb{Q}_{\nu}$-almost every $(\omega, x) \in \Omega \times M$, $\sigma_x^k$-almost every $v \in \mathrm{Gr}_k(T_x M)$
    \begin{equation}
        \lim_{t\to \infty}\frac{1}{t}\log r_t^k (\omega, x, v)= \lim_{t\to\infty}\frac{1}{t}\log\left\|\Bigwedge^k \Phi_t(\omega, x) v\right\| = \Lambda_1 +\cdots+ \Lambda_k .\label{rlim}
    \end{equation}
\end{corollary}

\begin{proof}

    The proof follows in the exact same way as the one of \cite[Corollary 2.1]{Baxendale1986TheDiffeomorphisms}, replacing the measure $\rho\otimes \mathbb P$ defined above \cite[Theorem 2.1]{Baxendale1986TheDiffeomorphisms} by the measure $\mathbb Q_\nu.$
    
    Equivalently, introducing the Borel measure $\mathcal{V}^k$ on $\Omega \times \mathrm{Gr}_k(M)$ such that
$$\mathcal{V}^k(A \times B) = \int_{A\times M}  \sigma_x^k\left( B\cap  \mathrm{Gr}_k(T_x M) 
\right) \mathbb Q_{\nu}(\d \omega, \d x), \ \forall \ A\times B \in \mathcal F \otimes \mathcal B (\mathrm{Gr}_k(M))$$

    we have
    \begin{equation*}
        \lim_{t\to \infty}\frac{1}{t}\log r_t^k (\omega, x, v) = \Lambda_1 +\cdots+ \Lambda_k \quad \mathcal{V}^k\text{-almost surely}.
    \end{equation*}
    This finishes the proof.
\end{proof}
\begin{remark}
    Observe that $r_t^k(\omega,x,v)$ was originally defined on (\ref{rands}) just for values of $v \in \Bigwedge^k_0(T_x M).$ However it is clear that if $v_1,v_2 \in \Bigwedge^k_0(T_x M)$ are linearly dependent then
    $$\lim_{t\to \infty}\frac{1}{t}\log r_t^k (\omega, x, v_1) = \lim_{t\to \infty}\frac{1}{t}\log r_t^k (\omega, x, v_2),  $$
    implying that the limit (\ref{rlim}) is well defined for $v \in 
    \bigwedge^k (T_xM)$.
\end{remark}

We can now prove Proposition \ref{thm:L1QConv}.

\begin{proof}[Proof of Proposition \ref{thm:L1QConv}] We use estimates from the proof of \cite[Theorem 3.3.3]{Arnold1998RandomSystems}. Observe that for every $(x,\omega)\in \Omega\times M$ and  $v \in \mathrm{Gr}_k(\mathbb R^d)$ we obtain\label{proof:estimates}
\begin{equation}
    \left|\log\left\|\Bigwedge^k \Phi_t v\right\|\right| \leq \max\left\{\left|\log\left\|\Bigwedge^k \Phi_t\right\|\right|, \left|\log\left\|\Bigwedge^k \Phi_t^{-1}\right\|\right|\right\}.
    \label{pontosdeexclamacao}
\end{equation}
But by \cite[Theorem 3.3.3, Proof of Part (B)(b)]{Arnold1998RandomSystems}, subbaditivity and $\Theta_t$-invariance of $\mathbb{Q}_{\nu}$, we observe
\begin{align*}
\frac{1}{t}\mathbb{E}^{\mathbb{Q}}_{\nu}\left[\left|\log\left\|\Bigwedge^k\Phi^{\pm1}_t\right\|\right|\right]\    &\leq\frac{1}{t} \sum_{n = 0}^{\lfloor t\rfloor} \mathbb{E}^{\mathbb{Q}}_{\nu}\left[\sup_{0\leq\ell\leq 1}\left|\log\left\|\Bigwedge^k \Phi^{\pm1}_\ell \circ \Theta_n\right\| \right|\right]\\
&\leq2\mathbb{E}^{\mathbb{Q}}_{\nu}\left[\sup_{0\leq\ell\leq 1}\left|\log\left\|\Bigwedge^k \Phi^{\pm 1}_\ell\right\| \right|\right]\\
&\leq 2k\left(\mathbb{E}^{\mathbb{Q}}_{\nu}\left[\sup_{0\leq\ell\leq 1}\log^{+}\left\| \Phi_\ell\right\|+\sup_{0\leq\ell\leq 1}\log^{+}\left\| \Phi^{-1}_\ell\right\|\right]\right)<\infty,
\end{align*}
where we have used the integrability conditions \eqref{eqn:IC+} and \eqref{eqn:IC-}. Now, combining the above equation~\eqref{pontosdeexclamacao}, Theorem~\ref{thm:IntroErgMeas}, Corollary \ref{BaxendaleResult} and \cite[Part Three-Theorem 4.18]{Royden1988RealAnalysis}, we obtain the existence of a set $\widetilde{G} \subset \mathrm{Gr}_k(M)$ such that
$\rho^k (\widetilde{G}) = 1$ and
$$\lim_{t\to\infty}\mathbb{E}^{\mathbb{Q}}_{\nu}\left[ \left|\lambda^{(k)} - \frac{1}{t}\log\left\|\Bigwedge^k \Phi_tv\right\|\right| \right] = 0\qquad \text{for all }(x,v)\in \widetilde{G}.$$
This finishes the proof of this proposition.
\end{proof}

\newpage
\section{Convergence of finite-time Lyapunov exponents}
\label{sec:3.5}

We start this section with the proof of Theorem \ref{thm:L1Conv}.

\begin{proof}[Proof of Theorem \ref{thm:L1Conv}]

We divide the proof into five steps.

\begin{step}{1} We show that if there exists $p \in (1, \infty]$ such that
\begin{equation*}
\sup_{t\geq 0}  \mathbb E_x\left[|\Gamma_t|^{p}\mid \tau > t\right] <\infty,
\end{equation*}
then
$$ \lim_{a\to \infty}\sup_{t\geq 0}  \mathbb E_x\left[\mathbbm 1_{|\Gamma_t|> a} |\Gamma_t|\mid \tau > t\right]=0.$$

It should be noted that the above step can be viewed as a conditioned version of the de la Vallée Poussin principle \cite{Klenke2020ProbabilityTheory}. 
\end{step}
By a direct computation, we obtain
$$|\Gamma_t|^{p} \geq a^{p-1}\mathbbm 1_{|\Gamma_t|> a} |\Gamma_t|$$
and, thus,
$$ \lim_{a\to \infty}\sup_{t\geq 0}  \mathbb E_x\left[\mathbbm 1_{|\Gamma_t|> a} |\Gamma_t|\mid \tau > t\right] \leq \lim_{a\to \infty}a^{1-p}\sup_{t\geq 0}  \mathbb E_x\left[|\Gamma_t|^{p}\mid \tau > t\right]=0.$$

\begin{step}{2}
We show that if $\{\Gamma_t\}_{t\geq 0}$ is a family of random variables fulfilling the assumptions in $(ii)$ and \begin{align}
\sup_{t\geq 0}  \left\| \mathbbm{1}_{\left\{\tau>t\right\}}\Gamma_t\right\|_{L^\infty (\Omega\times M, \mathbb P_x)} < \infty,    \label{eq:dennis}
\end{align}
then
$$\lim_{t\to\infty} \mathbb E_x\left[\left|\Gamma_t - \Gamma^\star\right| \mid \tau >t\right] = 0. $$
\end{step}

Since $|\Gamma_t - \Gamma^{\star}|$ is a $\mathcal G_t$-random variable,
Proposition \ref{marathon} implies
\begin{align*}
\lim_{t\to\infty}\frac{e^{\beta  t}}{\eta(x)}\mathbb E_x\left[|\Gamma_t - \Gamma^{\star}| \cdot \eta \circ \varphi_t\right]& =  \lim_{t\to\infty}\mathbb Q_x\left[|\Gamma_t - \Gamma^{\star}|\right] = 0.
\end{align*}

Since for every $x\in M$
$$ \lim_{t\to\infty}\frac{e^{\beta  t} }{\eta(x)} \mathbb P_x [\tau > t] = 1,$$
we obtain
\begin{align}
\mathbb E_x\left[|\Gamma_t - \Gamma^{\star}| \cdot \eta \circ \varphi_t \mid \tau >t \right]&= 0. \label{E41}
\end{align}

Moreover, given $\delta >0$ and $t\in\mathbb T,$ 
\begin{align}
    \delta\cdot \mathbbm 1_{\{\eta > \delta\} }\circ \varphi_t\cdot \mathbbm 1_{\{\tau>t\}} \leq \eta \circ \varphi_{t}\cdot \mathbbm 1_{\{\tau >t\}}= \eta\circ \varphi_{t}.\label{E51}
\end{align}

Combining equations \eqref{E41}--\eqref{E51}  yields
\begin{align*}
    &\lim_{t\to\infty}\mathbb E_x \left[ |\Gamma_t - \Gamma^{\star}|\cdot \mathbbm 1_{\{\eta > \delta\} }\circ \varphi_t \mid \tau >t\right] \leq \frac{1}{\delta}\lim_{t\to\infty}\mathbb E_x\left [  \left|\Gamma_t - \Gamma^{\star}\right|\cdot \eta\circ \varphi_t \mid \tau >t\right] = 0. 
\end{align*}

Observe that, defining $$K =|\Gamma^{\star}|+\sup_{t\geq 0} \left\|\mathbbm 1_{\{\tau > t\}} \Gamma_t\right\|_{L^\infty(\Omega\times M, \mathbb P_x)}<\infty,$$ we obtain
\begin{align*}
\mathbb E_x \left[ \left|\Gamma_t - \Gamma^{\star}\right|\cdot \mathbbm 1_{\{\eta \leq \delta\} }\circ \varphi_t  \mid\tau >t\right]\leq&  K   \mathbb E_x[\mathbbm 1_{\{\eta \leq \delta\} }\circ \varphi_t  \mid \tau >t].
\end{align*}

The above equation implies that
$$ \lim_{t\to\infty} \mathbb E_x \left[ \left|\Gamma_t - \Gamma^{\star}\right|\cdot \mathbbm 1_{\{\eta \leq \delta\} }\circ \varphi_t  \mid\tau >t\right] = K \mu(\{\eta\leq \delta\}). $$

This implies that for every $\delta>0,$
\begin{align*}
    0&\leq \limsup_{t\to\infty}\mathbb E_x\left[ \left|\Gamma_t - \Gamma^{\star}\right|\mid\tau>t\right]\\
    &=  \limsup_{t\to\infty}\left\{\mathbb E_x\left[ \left|\Gamma_t - \Gamma^{\star}\right|\mathbbm 1_{\{\eta >\delta\}}\circ \varphi_t \mid\tau>t\right]
    + \mathbb E_x\left[\left|\Gamma_t - \Gamma^{\star}\right| \mathbbm 1_{\{\eta \leq \delta\}}\circ \varphi_t \mid\tau>t\right]\right\} \\
    &\leq K\mu(\{\eta \leq \delta\}). 
\end{align*}
Since $\delta$ is arbitrary small and
$$\mu(\{\eta = 0\}) = 0, $$
we obtain
$$\lim_{t\to\infty}\mathbb E_x\left[\left|\Gamma_t - \Gamma^{\star}\right|\mid \tau>t\right] = 0.$$

\begin{step}{3}
Let $\{\Gamma_t\}_{t\geq 0}$ be a sequence of random variables fulfilling the assumptions of $(ii)$. Define $\Gamma_t^{(a)} := \mathbbm 1_{|\Gamma_t| \leq a} \Gamma_t.$ We show that for each $a\geq 2|\Gamma^\star|$ we have
\begin{equation}\label{gammaAConv}
    \lim_{t\to \infty} \mathbb E_x \left[\left|\Gamma^{(a)}_t - \Gamma^{\star}\right| \mid \tau >t\right] = 0.
\end{equation}
\end{step}
First observe that the restriction $a\geq 2|\Gamma^\star|$ implies $|\Gamma^{(a)}_t - \Gamma^{\star}| \leq |\Gamma_t - \Gamma^{\star}|$ and thus by (\ref{normal}) we also have
$$\lim_{t\to \infty} \mathbb{E}_{x}^{\mathbb Q}\left[\left|\Gamma^{(a)}_t - \Gamma^{\star}\right|\right] =0.$$
Furthermore, since $\left|\Gamma^{(a)}_t - \Gamma^{\star}\right| \leq a + |\Gamma^\star|$, equation \ref{eq:dennis} is fulfilled and from Step 2 we obtain
$$\lim_{t\to \infty} \mathbb E_x \left[\left|\Gamma_t^{(a)} - \Gamma^{\star}\right|\Big|~\tau >t\right] = 0.$$

\begin{step}{4}
    We show $(ii).$
\end{step}
    Observe that for each $a \geq 2|\Gamma^\star|$ we have
    \begin{align*}
        \limsup_{t\to \infty} \mathbb E_x \big[|\Gamma_t - \Gamma^{\star}|\mid \tau >t\big]&\leq \limsup_{t\to \infty} \left(\mathbb E_x \left[\left|\Gamma_t^{(a)} - \Gamma^{\star}\right|\Big|~  \tau >t\right] + \mathbb E_x \left[\left|\Gamma_t - \Gamma_t^{(a)}\right|\Big|~  \tau >t\right]\right)\\
        &= \limsup_{t\to \infty} \mathbb E_x \big[\mathbbm 1_{|\Gamma_t|> a} |\Gamma_t|\mid \tau >t\big] \\
        &\leq \sup_{t\geq 0}  \mathbb E_x\big[\mathbbm 1_{|\Gamma_t|> a} |\Gamma_t|\mid \tau > t\big].
    \end{align*}
    Considering $a \to \infty$ yields 
    $$ \lim_{t\to \infty} \mathbb E_x \big[|\Gamma_t - \Gamma^{\star}|\mid \tau >t\big] = 0,$$
    completing the proof.

\begin{step}{5}
    We show $(i).$
\end{step}
    Given $\varepsilon>0,$ consider the family of random variables
$$\left\{\widetilde \Gamma_t = \mathbbm 1_{\left\{|\Gamma_t - \Gamma^*|>\varepsilon\right\}} \right\}_{t\geq 0},$$
and define $\widetilde\Gamma^* = 0.$ From Step 4, we obtain
$$\lim_{n\to \infty}\mathbb P_x \left[|\Gamma_t - \Gamma^*|>\varepsilon \mid \tau >t\right] =  \lim_{n\to \infty}\mathbb E_x \left[\left| \widetilde{\Gamma}_t - \widetilde{\Gamma}^* \right|\Big|~  \tau >t\right] =0, $$
proving the theorem.
\end{proof}

In the following, we use Theorem \ref{thm:L1Conv} to prove Corollaries \ref{T2.11} and \ref{T2.12}.

\begin{proof}[Proof of Corollary \ref{T2.11}]

Let $\varepsilon > 0$. From Corollary \ref{BaxendaleResult}, we obtain that for $\rho^k$-almost every~${(x,v)\in \mathrm{Gr}_k(M)}$, 
\begin{align*}
    \lim_{t\to\infty}\mathbb Q_x\left[ \left\{\left|\lambda^{(k)} - \frac{1}{t} \log\left\|\Bigwedge^k \Phi_t v\right\|\right|> \varepsilon \right\}\right] = 0. 
\end{align*}

Let us define $\Gamma^* = \lambda^{(k)}$ and for every $t\geq 0$ the $\mathcal G_t$-measurable random variable
$$ \Gamma_t = \frac{1}{t} \log\left\|\Bigwedge^k \Phi_t v\right\|.$$
Applying Theorem \ref{thm:L1Conv} to the family $\{\Gamma_t\}_{t\geq 0}$, we immediately obtain the desired result and similarly to show the convergence \eqref{eqn:CondProbNorm}.
\end{proof}

\begin{proof}[Proof of Corollary \ref{T2.12}]

From Proposition~\ref{thm:L1QConv}, we immediately get that there exists a subset ${\widetilde{G}\subset \mathrm{Gr}_k(M)}$ of full $\rho^k$-measure such that for all $(x,v) \in \widetilde{G}$
$$\lim_{t\to\infty}\mathbb{E}_{x}^{\mathbb Q}\left[ \left|\lambda^{(k)} - \frac{1}{t}\log r^k_t(v)\right|  \right] = 0.$$

Therefore, taking $\Gamma^{\star} = \lambda^{(k)}$ and the $\mathcal{G}_t$-measurable random variable
$$
\Gamma_t = \frac{1}{t}\log\left\|\Bigwedge^k \Phi_t(\cdot, x) v\right\| ,
$$
in combination with \cite[Theorem 3.3.3, Proof of Part (B)(b)]{Arnold1998RandomSystems} and (\ref{Qnu}), yields
$$\sup_{t\geq 0} \|  \Gamma_t \|_{L^p (\Omega\times M, \mathbb P_x(\cdot \mid \tau >t))} <\infty,\quad \text{for}\ \nu\text{-almost every }x\in M. $$

Applying Theorem~\ref{thm:L1Conv} immediately yields the desired result and similarly for \eqref{C5}. This completes the proof of the theorem.

\end{proof}

\newpage

\section{Application to Stochastic flows}
\label{sec:4}

\label{PT2.10}
Let us recall that we work on the probability space
\begin{align*}
    (\Omega, (\mathcal{F}_t)_{t\geq 0}, \mathcal{F},\mathbb P) = \left(\mathcal C_0(\R_+,\R^m),\left(\sigma(\pi_s,0\leq s\leq t)\right)_{t\geq 0}, \sigma(\pi_s,s\geq 0),\mathbb P\right),
\end{align*}

as in Eq. \eqref{l1} where $\Pb$ is the Wiener measure.

Let us assume that equation $(\ref{SDE})$ satisfies Hypothesis \ref{(B)}. By \cite[Theorem 1.2.9.]{Hsu2002StochasticManifolds}, then for each initial condition $x\in M$, there exists a unique solution ${\varphi_t(\cdot, x): \mathbb R_+\times \Omega \to E_M = M\sqcup \{\partial\}}$ of the stochastic differential equation
\begin{align}
    \d X_t  &= V_0 (X_t) \d t + \sum_{i=1}^{m} V_i(X_t) \circ \d W^i_t,\ X_0=x, \label{SDE1}
\end{align} 
on $M$ up until explosion, i.e.~defining the stopping time
$$\tau(\omega,x) = \inf \{t\geq 0 : \varphi_t(\omega, x) \not\in M\}. $$
The stochastic process $\varphi_t(\cdot,x)$ fulfils the following conditions
\begin{enumerate}
    \item[(i)] $\varphi_0(\cdot ,x) =x;$
    \item[(ii)] $\displaystyle \varphi_t(\cdot, x) = \int_0^t V_0(\varphi_s(\cdot, x))\d s + \sum_{i=1}^m \int_0^t V_i(\varphi(\cdot, x)) \d s, \ \forall \ 0\leq t\leq \tau(\cdot,x);$
    \item[(iii)] for every $t\geq 0$, $\varphi_t(\cdot, x)$ is $\mathcal F_t$-measurable;
    \item[(iv)] $\varphi_t(\cdot,x) = \partial,$ for every $\tau(\cdot, x) \leq t;$
    \item[(v)] for every $\omega\in \Omega,$ the path $t \mapsto \varphi_t(\omega, x)$ is continuous on $E_M,$
    \item[(vi)] if $Y_t:\Omega \to E_M$ is another stochastic process satisfying conditions (i)--(v) then
    $$\mathbb P[Y_t = \varphi_t(\cdot, x), \ \forall \ t\geq 0] = 1. $$
\end{enumerate}

Thus, this defines a stochastic flow which we denote by the same symbol $\varphi$
\begin{align*}
    \varphi: \mathbb{R}_{+}\times \Omega \times E_M &\to E_M\\
    (t,\omega,x)&\mapsto \varphi_t(\omega, x).
\end{align*}
From \cite[Chapter 2.3]{Arnold1998RandomSystems} and \cite[Proposition 2.5]{Ling2021TheSDEs}, we may assume without loss of generality that $\varphi$ forms a \textit{perfect} cocycle (by taking a modification, see Remark \ref{rmk:modification}), i.e.
$$\varphi_{t+s} (\omega,x) = \varphi_t ( \theta_s \omega,  \varphi_s(\omega,x)),\ \text{for every }s,t\geq 0,$$
where
\begin{align*}
    \theta_t: \Omega &\to \Omega\\
    \omega &\mapsto \omega (t + \cdot) - \omega(t).
\end{align*}

With the above notation, we say that $(\theta,\varphi)$ is the random dynamical system induced by the stochastic differential equation \eqref{SDE}.

From \cite[Theorem 2.3.32]{Arnold1998RandomSystems}, we obtain that the linearised flow $\Phi_t(\omega,x):=\mathrm D \varphi_t (\omega,x)$ solves the stochastic differential equation
$$\begin{cases}
    \displaystyle \d \Phi_t = \mathrm{D}V_0(\varphi_t ) \Phi_t \d t + \sum_{i=1}^m \mathrm{D}V_i(\varphi_t) \Phi_t \circ \d W_t^i, \qquad \forall \ 0\leq t\leq \tau, \\
    \Phi_0 = \mathrm{Id}.
\end{cases}$$

Fixing $k\leq d$ and $v = v_1 \wedge\cdots \wedge v_k\in \mathrm{Gr}_k(T_x M) \simeq \mathrm{Gr}_k(\R^d)$, recall that we denote
\begin{align*}
    r^k_t(\omega, x, v) := \left\|\Bigwedge^k\Phi_t(\omega, x) v\right\|\qquad \text{and} \qquad s_t^k(\omega, x, v) = \frac{\Bigwedge^k\Phi_t(\omega, x) v}{\left\|\Bigwedge^k\Phi_t(\omega, x) v\right\|}\in \mathrm{Gr}_k(T_x M)
\end{align*}

From \cite[Theorem 3.1]{Baxendale1986TheDiffeomorphisms}, there exist continuous (and hence bounded) functions ${\psi^k: \mathrm{Gr}_k(M)\to\mathbb R}$ and $\phi^k_i:\mathrm{Gr}_k(M)\to\mathbb R$ for $i\in\mathbb \{1,\dots,m\}$ such that
\begin{align}
    \d \left(\log r_t^k\right) = \psi^k\left(\varphi_t, s_t^k\right) \d t+ \sum_{i=1}^{m}\phi_i^k(\varphi_t, s^k_t) \d W_t^i, \ \forall \ 0\leq t\leq \tau(\cdot ,x).\label{f1} 
\end{align}
The Formulae for these functions were derived in \cite{Baxendale1986TheDiffeomorphisms} and are given by
$$\psi^k(x,s) := \mathrm{tr}(V'_0(x)P_s) + \sum_{i = 1}^m\left\{\mathrm{tr}(V_i'(x)V_i(x)P_s) -\mathrm{tr}(V_i'(x)P_s)^2+\mathrm{tr}(V_i'(x)^{*}(I-P_s)V_i'(x)P_s)\right\}$$
and
$$\phi_i^k(x,s) := \mathrm{tr}(V_i'(x)P_s)$$
where $P_s$ denotes the projection onto the subspace $s\in \mathrm{Gr}_k(\R^d)$.
We can now prove Theorem \ref{T2.10}.

\begin{proof}[Proof of Theorem \ref{T2.10}]
{\revision From \cite[Corollary 1.9]{Benaim2021DegenerateDomain} and the computation done in \cite[Appendix]{oçafrain2023central}},
it is clear that the random dynamical system $(\theta,\varphi)$ associated to $(\ref{SDE2})$ satisfies Hypothesis $\mathrm{\ref{(H)}}$. In the remainder of this proof, in the interest of readability, we fix $k\leq d$ and drop this superscript.

It suffices for us to show that the integrability conditions \eqref{eqn:IC+} and \eqref{eqn:IC-} are fulfilled allowing us to apply Theorem \ref{thm:fullFK}. Moreover, for all $T>0$ and for every $k\in \{1,\ldots,d\}$,
    \begin{align}
        \sup_{0\leq t\leq T} \log \|\Bigwedge^k\Phi_t\|,\  \displaystyle\sup_{0\leq t\leq T} \log \|\Bigwedge^k \Phi_t^{-1}\| \in L^1(\Omega\times M, \mathcal F\otimes \mathcal B(M), \mathbb Q_\nu),\label{step0eq}
    \end{align}

We follow ideas of \cite[Remark 6.2.12]{Arnold1998RandomSystems}. First note that for $A \in GL(d, \mathbb{R})$
$$\frac{\left\|\left(\Bigwedge^k A\right) (e_{i_1}\wedge \cdots \wedge e_{i_k})\right\|}{\left\|e_{i_1}\wedge \cdots \wedge e_{i_k}\right\|}\leq\left\|\Bigwedge^k A\right\| \leq \begin{pmatrix}d\\k \end{pmatrix}\max_{1\leq i_1\leq \cdots i_k \leq d} \frac{\left\|\left(\Bigwedge^k A\right) (e_{i_1}\wedge \cdots \wedge e_{i_k})\right\|}{\left\|e_{i_1}\wedge \cdots \wedge e_{i_k}\right\|}.$$

So, for the first integrability condition, it suffices for us to prove that for $v \in \mathrm{Gr}_k(\R^d)$,

\begin{equation}
    \mathbb{E}^{\mathbb{Q}}_{\nu}\left[\sup_{0\leq t\leq T}|\log r_t(\cdot, \cdot, v)|\right]<\infty\label{eqn:ICMETvec}
\end{equation}
But for $t \leq T$,
\begin{equation}
|\log r_t|=\left|\int_{0}^t \psi(\varphi_\ell, s_\ell)\d \ell +\sum_{i=1}^d\int_{0}^t\phi_i(\varphi_\ell, s_\ell) \d W^i_\ell\right|\\\leq t\|\psi\|_{\infty} +\sum_{i=1}^d\left|\int_{0}^t\phi_i(\varphi_\ell, s_\ell) \d W^i_\ell\right|\label{eqn:ICsdeproof}
\end{equation}

Now for $ i \in \{ 1, \ldots, d\}$
\begin{align*}
    \mathbb{E}^{\mathbb{Q}}_{\nu}\left[\sup_{0\leq t\leq T}\left|\int_0^t \phi_i(\varphi_\ell, s_\ell)\d W^{i}_\ell\right|\right] &= \int_M \mathbb{E}_{x}^{\mathbb{Q}}\left[\sup_{0\leq t\leq T}\left|\int_0^t \phi_i(\varphi_\ell, s_\ell)\d W^{i}_\ell\right|\right]\nu(\d x)\\
    &\leq \int_{M}\frac{e^{\beta  T}}{\eta(x)}\mathbb{E}_x\left[\eta \circ\varphi_t\sup_{0\leq t\leq T} \left|\int_0^t \phi_i(\varphi_\ell, s_\ell)\d W^{i}_\ell\right| \mathbbm{1}_{\{\tau >T\}}\right]\eta(x)\mu(\d x)\\
    &\leq\|\eta\|_{\infty} e^{\beta  T}\int_{M}\mathbb{E}_x\left[\sup_{0\leq t\leq T} \left|\int_0^{t} \phi_i(\varphi_\ell, s_\ell)\d W^{i}_\ell\right| \mathbbm{1}_{\{\tau >T\}}\right]\mu(\d x),
\end{align*}
since $\sup_{0\leq t\leq T}\left|\int_0^t \phi_i(\varphi_\ell, s_\ell)\d W^{i}_\ell\right|^2$ is $\mathcal{F}_T$-measurable. Now,
\begin{align*}    
    \mathbb{E}_x\left[\sup_{0\leq t\leq T} \left|\int_0^t \phi_i(\varphi_\ell, s_\ell)\d W^{i}_\ell\right| \mathbbm{1}_{\{\tau >T\}}\right]&\leq\mathbb{E}_x\left[\sup_{0\leq t\leq T} \left|\int_0^{t\wedge \tau} \phi_i(\varphi_\ell, s_\ell)\d W^{i}_\ell\right|\right]\\
   &\leq
     \left(\mathbb{E}_x\left[\sup_{0\leq t\leq T} \left|\int_0^{t\wedge \tau} \phi_i(\varphi_\ell, s_\ell)\d W^{i}_\ell\right|^2\right]\right)^{1/2},
\end{align*}
where we use H\"{o}lder's\ inequality in the last step.
Now observe that $\int_0^{t\wedge \tau} \phi_i(\varphi_\ell, s_\ell)\d W^{i}_\ell$ is a stopped martingale. Therefore, by the Burkholder--Davis--Gundy inequality \cite[Chapter IV, Corollary 4.2]{Revuz1999ContinuousMotion},

\begin{equation}
\left(\mathbb{E}_x\left[\sup_{0\leq t\leq T} \left|\int_0^{t\wedge \tau} \phi_i(\varphi_\ell, s_\ell)\d W^{i}_\ell\right|^2\right]\right)^{1/2}\leq\left(4 \mathbb{E}_x\left[\int_0^{T\wedge \tau}|\phi_i(\varphi_\ell, s_\ell)|^2\d \ell\right]\right)^{1/2}\leq 2 T^{1/2} \|\phi_i\|_{\infty}.
\end{equation}

Since this bound is uniform over all $x \in M$, this finishes the proof that \eqref{eqn:ICMETvec} holds uniformly over all $v \in \mathrm{Gr}_k(\R^d)$. This proves the first integrability condition of equation \eqref{step0eq}.

Now, observe that $\Phi_t^{*-1}$ solves the stochastic differential equation

$$\begin{cases}
    \displaystyle \d \Phi^{*-1}_t = -\mathrm{D}^{*}V_0(\varphi_t ) \Phi^{*-1}_t \d t - \sum_{i=1}^m \mathrm{D}^{*}V_i(\varphi_t) \Phi^{*-1}_t \circ \d W_t^i, \qquad \forall \ 0\leq t\leq \tau, \\
    \Phi^{*-1}_0 = \mathrm{Id}.
\end{cases}$$

Since $\left\|\Phi^{-1}_t\right\| = \left\|\Phi_t^{*-1}\right\|$, applying the same reasoning as above proves the second integrability condition of equation \eqref{step0eq}.

These immediately imply integrability conditions \eqref{eqn:IC+} and \eqref{eqn:IC-}. Thus, Theorems \ref{thm:fullFK} and \ref{thm:fullMET} and Corollary \ref{T2.11}  hold and yield the existence of conditioned Lyapunov exponents ${\Lambda_1 \geq \cdots \geq \Lambda_d> -\infty}$.
\begin{equation*}
    \lim_{t\to\infty}\mathbb{E}^{\mathbb{Q}}_{\nu}\left[\left|\lambda^{(k)} - \frac{1}{t}\log\left\|\Bigwedge^k \Phi_t\right\|\right|\right] = 0 
\end{equation*}
Furthermore, by an application of the subadditive ergodic theorem (similarly to \cite[Chapter 3]{Viana2014LecturesExponents}), we have
\begin{equation*}
    \lim_{t\to\infty}\mathbb{E}^{\mathbb{Q}}_{\nu}\left[\left| (\Lambda_{d-k+1} +\cdots +\Lambda_d) + \frac{1}{t}\log\left\|\Bigwedge^k \Phi^{*-1}_t\right\|\right|\right] = 0 
\end{equation*}

Corollary \ref{T2.11} proves the desired convergence in conditional probability: for every $\varepsilon> 0$ and $\rho^k$-almost every $(x,v) \in \mathrm{Gr}_k(M),$
\begin{align*}
    \lim_{t\to\infty } \mathbb P_x \left[ \left\{ \left|\lambda^{(k)} - \frac{1}{t}\log \|\Bigwedge^k \Phi_t v\|\right| >\varepsilon \right\}\bigg|~  \tau > t \right]  = 0,
\end{align*}
which finishes the proof.
\end{proof}

Furthermore, we wish to generalise the Furstenberg–Khasminskii formula given by \cite{Engel2019ConditionedSystems} to compute the top Lyapunov exponent and provide an equivalent for lower exponents. This is, once more, achieved by following ideas of Baxendale \cite{Baxendale1986TheDiffeomorphisms} below. A remarkable aspect of these formulae is the apparent impossibility to derive the multiplicative noise case without the use of the $Q$-process. This additionally demonstrates the usefulness of the $Q$-process to study conditioned finite-time dynamics.

\begin{proposition}
    For $k\leq d$, if there exist unique quasi-stationary and unique quasi-ergodic distributions $\mu^k$ and $\nu^k$ on $\mathrm{Gr}_k(M)$ for the process $(\varphi_t, s^k_t)$, then
    
    \begin{equation}
    \lambda^{(k)} = \Lambda_1 +\cdots +\Lambda_k  = \int_{\mathrm{Gr}_k(M)} \psi^k\d \nu^k + \sum_{i = 1}^m\sum_{j =1}^d\int_{\mathrm{Gr}_k(M)}\phi^k_i V^j_i \partial_j\eta\d \mu^k\label{eqn:sde-qed}
    \end{equation}
    and similarly for the last $k$ Lyapunov exponents. In particular,
    \begin{equation*}
    \lambda^{\pm} = \int_{\mathbf{P}^{d-1}(TM)} \psi^{\pm 1}\d \nu^{\pm 1} + \sum_{i = 1}^m\sum_{j =1}^d\int_{\mathbf{P}^{d-1}(TM)}\phi^{\pm 1}_i V^j_i\partial_j\eta\d \mu^{\pm 1}
    \end{equation*}
    where $\lambda^{+}, \lambda^{-}$ denote the extremal Lyapunov exponents $\Lambda_1$ and $\Lambda_d$ respectively and $\psi^{\pm1}, \phi_i^{\pm1}$ are given by \eqref{step0eq} and its analogue for the inverse linearised flow and similarly for $\nu^{\pm1}$ and $\mu^{\pm1}$.
    \label{cor:sde-qed}
\end{proposition}
 
\begin{proof}
Recall that for $\mathcal{V}^k$-almost every $(\omega, x, v)\in \Omega \times \mathrm{Gr}_k(M)$,
    
$$\lambda^{(k)} = \Lambda_1 +\cdots +\Lambda_k = \lim_{t\to\infty}\frac{1}{t}\log r^k_t(\omega, x, v)$$
    Also recall, formula \eqref{f1} below (where we again drop the superscript $k$)
\begin{align*}
    \d \left(\log r_t\right) = \psi\left(\varphi_t, s_t\right) \d t+ \sum_{i=1}^{m}\phi_i(\varphi_t, s_t) \d W_t^i, \ \forall \ 0\leq t\leq \tau(\cdot ,x).
\end{align*}

The time-average of the first term converges by Birkhoff's ergodic theorem
$$\lim_{t\to\infty}\mathbb{E}^{\mathbb{Q}}_{\nu}\left[\frac{1}{t}\int_0^t \psi(\varphi_\ell, s_\ell) \d \ell\right] = \int_{\mathrm{Gr}_k(M)} \psi \d \nu^k.$$

Now for the second term, by Girsanov's Theorem one can show the existence of a $h$-transform
$$\d W^i_t = \sum_{j}^d V^j_i(\varphi_t(\cdot, x))\partial_j\left(\log\eta(\varphi_t(\cdot, x))\right) \d t+\d B^i_t$$
where $\left(B^i_t\right)_{t\geq 0}$ is a $\mathbb{Q}_x$-standard Brownian motion. Now, on the one hand, by H\"{o}lder's inequality and  It\^{o} isometry, we obtain
\begin{align*}
    \left|\frac{1}{t}\mathbb{E}_{x}^{\mathbb{Q}}\left[\int_0^t \phi_i(\varphi_\ell, s_\ell) \d B^i_\ell\right]\right| &\leq \frac{1}{t}\left(\mathbb{E}_{x}^{\mathbb{Q}}\left[\left|\int_0^t \phi_i(\varphi_\ell, s_\ell) \d B^i_\ell\right|^2\right]\right)^{1/2}\\
    &\leq \frac{1}{t}\left(\mathbb{E}_{x}^{\mathbb{Q}}\left[\int_0^t \phi_i^2(\varphi_\ell, s_\ell) \d \ell\right]\right)^{1/2}\\
    &\leq \frac{\|\phi_i\|_{\infty}}{\sqrt{t}}\xrightarrow[t \to\infty]{} 0
\end{align*}

On the other hand,

\begin{align*}
    \mathbb{E}^{\mathbb{Q}}_{\nu}\left[\frac{1}{t}\int_0^t\phi_i(\varphi_\ell, s_\ell)  V_i^j(\varphi_\ell)\partial_j\left(\log\eta(\varphi_\ell)\right)\d \ell\right] \xrightarrow[t\to\infty]{\textrm{Birkhoff}}&\int_{\mathrm{Gr}_k(M)}\phi_i V_i^j\partial_j(\log\eta)\d \nu^k\\
    =&\int_{\mathrm{Gr}_k(M)} \phi_i V^j_i \partial_j\eta \d \mu^k.
\end{align*}
Yielding the desired result \eqref{eqn:sde-qed}. Note that, although it explodes near the boundary $\partial M$, the integrand $\phi_i V_i^j\partial_j(\log\eta) \in L^1(\nu^k)$ since $\d\nu^k = \eta \d\mu^k$.
\end{proof}

This proposition generalises the definition of conditioned Lyapunov exponents given in \cite{Engel2019ConditionedSystems} which treats the additive case with $k=1$. The Lyapunov exponents can then be computed recursively: $\Lambda_k = \lambda^{(k)}-\lambda^{(k-1)}$. This is particularly useful for numerical estimations of the conditioned Lyapunov exponents. Note that the process $(\varphi_t, s^k_t)$ is degenerate, making the uniqueness of its quasi-stationary and quasi-ergodic distributions unclear in general. Some criteria for the exponential convergence of this process to quasi-stationarity, such as the H\"{o}rmander condition, are discussed in \cite{Benaim2021DegenerateDomain}. 

A particular case of this proposition is the Liouville's formula below.

\begin{corollary}[Liouville's formula]
    Let $\mu$ and $\nu$ be the quasi-stationary and quasi-ergodic distributions of $\varphi$ on $M$, then
    $$\lambda^{(d)} =  \lim_{t\to\infty}\mathbb{E}^{\mathbb{Q}}_{\nu}\left[\frac{1}{t}\log\det(\Phi_t)\right]=\lim_{t\to\infty}\mathbb{E}^{\mathbb{Q}}_{\nu}\left[\frac{1}{t}\log\left\|\Bigwedge^d\Phi_t\right\|\right]=
    \int_{M} \psi^d \d \nu +\sum_{i=1}^m \sum_{j=1}^d \int_M \phi_i^dV_i^j \partial_j \eta \d \mu.$$\label{cor:Liuouville}
\end{corollary}

Finally, we give the corollary below as an application of Corollary \ref{T2.12} for absorbed diffusions with additive noise.
\begin{corollary}[Additive noise case]\label{cor:AddSDE}
    Let $(\theta, \varphi)$ be as in Theorem \ref{T2.10} and assume further that the vector fields $\left\{V_i\right\}_{i=1}^m$ are constants, i.e.~$(\theta, \varphi)$ is generated by a stochastic differential equation with additive noise
    \begin{align}
    \d X_t  = V_0 (X_t) \d t + \sum_{i=1}^{m} V_i \circ \d W_i^t,\ X_0=x,
\end{align} on $M$ up until explosion.
Then the convergence of the finite-time Lyapunov exponents occurs in conditional expectation in the sense that for all $k\leq d$, for $\rho^k$-almost every $(x,v) \in \mathrm{Gr}_k(M),$
\begin{align}
    \lim_{t\to\infty } \mathbb E_x \left[\left|\lambda^{(k)} - \frac{1}{t}\log \left\|\Bigwedge^k\Phi_t v\right\| \right| \bigg|~ \tau > t \right]  = 0. \label{HC1}
\end{align}
and

$$\lambda^{(k)} =\int_{\mathrm{Gr}_k(M)} \langle s, \widehat{\mathrm{D}V}_0^k(x)s\rangle \nu^k(\d x, \d s)$$
where for $A \in \mathrm{End}(\R^d)$, $\hat{A}^k \in \mathrm{End}(\Bigwedge^k\R^d)$ is defined on $\Bigwedge^k_0\R^d$ as
$$\hat{A}^k(v_1 \wedge \cdots \wedge v_k) := \sum_{i = 1}^k( v_1 \wedge \cdots \wedge A v_i \wedge \cdots \wedge  v_k).$$

\end{corollary}

\begin{proof}
    Observe that for all $k\leq d$ and for all $i \in \{1, \ldots, m\}$, $\phi_i^k = 0$ and ${\psi^k :(x, s)\mapsto  \langle s, \widehat{\mathrm{D}V}_0^k(x)s\rangle}$. Therefore, \eqref{eqn:ICsdeproof} directly implies that the integrability condition $\eqref{Qnu}$ is fulfilled and Corollary \ref{T2.12} yields the desired result.
\end{proof}

We mention that under additional assumptions, it is possible to characterise the fluctuations of the limit \eqref{HC1} via a central limit theorem for absorbed Markov processes (see for instance \cite[Theorem 1]{oçafrain2023central}).

\begin{appendix}
\section{Random Dynamical Systems}
\label{RDS}
In this appendix, we recall the definition of a random dynamical system. Let $\mathbb T$ be $\mathbb N_0$ or $\mathbb R_+.$ In the interest of clarity, our notations correspond to the ones of continuous time, e.g. sums over discrete time are denoted as integrals. 

\begin{definition}[Metric Dynamical System]
    $(\Omega, (\mathcal{F}_t)_{t\in \mathbb{T}}, \mathcal{F}, \Pb)$ be a filtered probability space. A family of mappings $\theta = \left\{ \theta_t : (\Omega, (\mathcal{F}_t)_{t\in \mathbb{T}}, \mathcal{F}, \Pb) \to (\Omega, (\mathcal{F}_t)_{t\in \mathbb{T}}, \mathcal{F}, \Pb)\right\}_{t\in \mathbb{T}}$ is said to be a metric dynamical system (or measure preserving DS) if it satisfies the following:

\begin{enumerate}
    \item $(\omega, t)\mapsto \theta_t\omega$ is $(\mathcal{F}\otimes\mathcal{B}(\mathbb{T})-\mathcal{F})$-measurable;
    \item $\theta_0 =\mathrm{id}_{\Omega};$
    \item Semiflow property: $\theta_{s+t} = \theta_s \circ \theta_t \ \textrm{for all $s, t\in \mathbb{T}$}$;
    \item $\Pb$ is $\theta_t$-invariant for all $t\in \mathbb{T}$, i.e.~$\left(\theta_{t}\right)_*\Pb = \Pb$ where $\left(\theta_{t}\right)_*\Pb (A) = \Pb(\theta_{t}^{-1}(A)),$ for all $A\in\mathcal F$;
    \item $\theta$ is said to be a filtered DS if $\theta_{t}^{-1}\mathcal{F}_s \subset \mathcal{F}_{s+t}$ for all $s,t \in \mathbb{T};$
    \item Furthermore, $\theta$ is said to be \textit{ergodic} if for all $t\in \mathbb{T}$, $\theta_t$-invariant sets have measure $0$ or $1$, i.e.~for all $A \in \mathcal{F}$, $\theta_{t}^{-1}A = A$ implies $\Pb(A) \in \{0,1\}.$
\end{enumerate}

When the context is clear, the quadruplet $(\Omega, \mathcal{F}, \Pb, (\theta_t)_{t\in \mathbb{T}})$ denotes a metric dynamical system $\theta$. If $\theta$ is a filtered DS, $(\Omega, (\mathcal{F}_t)_{t\in \mathbb{T}}, \mathcal{F}, \Pb, (\theta_t)_{t\in \mathbb{T}})$ might be referred to as the \textit{noise space}.
\end{definition}

We may impose the following additional condition on our noise space.

\begin{definition}[Memoryless Noise Space]
    A noise space $(\Omega, (\mathcal{F}_t)_{t\in \mathbb{T}}, \mathcal{F}, \Pb, (\theta_t)_{t\in \mathbb{T}})$ is said to be memoryless if for any $s, t \in \mathbb{T}$, $\mathcal{F}_s$ and $\theta_{s}^{-1}\mathcal{F}_t$ are independent under $\Pb$.\label{def:memoryless}
\end{definition}

We also recall the definition of a random dynamical system.

\begin{definition}
    A random dynamical system on a measurable state space $(X, \mathcal{B})$ over a metric dynamical system $(\Omega, \mathcal{F}, \Pb, (\theta_t)_{t\in \mathbb{T}})$ is a mapping
    
    \begin{align*}
        \varphi : \mathbb{T} \times \Omega \times X &\to X\\
        (t, \omega, x) &\mapsto \varphi_t(\omega, x)
    \end{align*}
   which satisfies the following properties
   \begin{enumerate}
       \item Measurability: $\varphi$ is a $(\mathcal{B}(\mathbb{R})\otimes \mathcal{F} \otimes \mathcal{B}-\mathcal{B})$-measurable mapping; and
       \item Cocycle property: $\varphi$ forms a \emph{perfect} cocycle over $\theta$, i.e.~$\omega \in \Omega$
        \begin{enumerate}
            \item $\varphi_0(\omega, \cdot)=\mathrm{id}_X$
            \item $\varphi_{t+s}(\omega, x)=\varphi_t(\theta_s\omega, \varphi_s(\omega, x))$ for all $s,t\in \mathbb{T}$ and for all $x\in X$.
        \end{enumerate}
    \end{enumerate}
  When the context is clear, we denote the random dynamical system $\varphi$ over the metric dynamical system $(\Omega,\mathcal F, \Pb, (\theta_t)_{t\in\mathbb T})$ simply by $(\theta, \varphi).$ When $(\Omega, (\mathcal{F}_t)_{t\in \mathbb{T}}, \mathcal{F}, \Pb, (\theta_t)_{t\in \mathbb{T}})$ is a memoryless noise space and $\varphi_t$ is $\mathcal F_t \otimes \mathcal B$-measurable for every $t\in\mathbb T, $  we say that $(\theta,\varphi)$ is a memoryless random dynamical system.\label{def:rds}
\end{definition}

\begin{definition}[RDS with absorption]
    Let $X$ be a topological state space that can be decomposed as $X =M\sqcup \{\partial\}$ where $M \subset X$ and $\{\partial \}$ denotes a so-called ``cemetery" or ``coffin" state. A measurable RDS $(\theta, \varphi)$ over a metric dynamical system $(\Omega, \{\mathcal F_t\}_{t\in\mathbb T} ,\mathcal{F}, \Pb, (\theta_t)_{t\in \mathbb{T}})$ is said to form a random dynamical system with absorption $(\theta, \varphi)$ on $X =M\sqcup \{\partial\}$ if for all $\omega\in \Omega$, $x \in X$, $\varphi_s(\omega, x)= \partial$ implies $\varphi_t(\omega, x) = \partial$ for all $t\geq s$. This justifies the definition of the following stopping time for each $x \in M$
    $$\tau(\cdot , x) = \inf\left\{t\geq 0 : \varphi_t(\cdot, x) = \partial\right\} $$
    In this context, a measurable RDS $(\theta, \varphi)$ is said to be
    \begin{itemize}
        \item continuous if for all $\omega\in \Omega$, the mappings
        \begin{align*}
            \varphi_{\cdot}(\omega, \cdot): \left\{(t, x) \in \mathbb{T}\times X \mid \tau(\omega, x)>t\right\} &\to M\\
            (t,x)&\mapsto\varphi_t(\omega, x)
        \end{align*}
        are continuous.
        \item If $X$ is furthermore endowed with a smooth structure, i.e.~if it is a manifold, then $(\theta, \varphi)$ is said to be of class $\mathcal C^k$ ($1\leq k\leq \infty$), if for all $t\in T$ and $\omega \in \Omega$, the mappings
        \begin{align*}
            \varphi_t(\omega)=\varphi_t(\omega, \cdot):  \left\{\tau(\omega, \cdot)>t\right\} \subset M&\to M\\
            x &\mapsto\varphi_t( \omega, x)
        \end{align*}
        are $k$-times differentiable (in the sense of \cite[Page 645]{Lee2012IntroductionManifolds}) where $\left\{\tau(\omega, \cdot)>t\right\} = \left\{ x \mid \tau(\omega, \cdot)>t\right\}.$
    \end{itemize}

\end{definition}

\section{Proof of Proposition \ref{prop:A}}
\label{Appendix:B}
In this section, we give a proof of proposition \ref{prop:A}. The proof below is based on the techniques developed in \cite{Champagnat2017UniformQ-process}, where similar results were proven assuming that the function $C(x)$ (in $\mathrm{(H2)}$ of Hypothesis $\mathrm{\ref{(H)}})$ is constant. 

\begin{proof}[Proof of Proposition \ref{prop:A}] 
The proof is done assuming $\mathbb{T}=\mathbb R_+$. If $\mathbb T=\mathbb N_{0}$, the same proof  holds with minor adaptations.  We divide the proof into four steps.

\begin{step}{1}\label{step1}
We show that for every non-negative measurable bounded function $g:M\to \mathbb R_+,$
$$\lim_{t\to\infty}e^{\beta  t}\mathcal P^{t}(g)(x) = \eta(x)\int_{M} g\  \d \mu, \ \text{for every} \ x\in M. $$
\end{step}

Since $g$ is non-negative and bounded, from $\mathrm{(H2)}$ we obtain that for every $t\geq 0,$
$$\left|\frac{\mathcal P^t (g)(x)}{\mathcal P^t (x,M)}  - \int_M g \ \d \mu\right|\leq \|g\|_{\infty}C(x)e^{-\alpha t}, $$
where $\|g\|_{\infty} := \sup_{x\in M} |g(x)|.$ Therefore, for every $t\geq 0,$
$$\left|e^{\beta  t}\mathcal P^t (g)(x) -   e^{\beta  t}\mathcal P^t (x,M) \int_M g \ \d \mu\right|\leq \|g\|_{\infty}C(x)e^{-\alpha t} e^{\beta  t}\mathcal P^t (x,M). $$
Since $$\lim_{t\to \infty}e^{\beta  t}\mathcal P^t(x,M) = \eta(x).$$
We obtain that
$$\lim_{t\to\infty}e^{\beta  t}\mathcal P^t (g)(x) =   \eta(x) \int_M g\  \d \mu.$$
This proves Step 1.

\begin{step}{2}
We prove $(i).$
\end{step}
Let $f:M\to\mathbb R$ be a non-negative measure function.  Let us consider the function
$$h_u(x) =  \min\left\{\eta(x),\inf_{r\geq u} \left\{e^{\beta  r}\mathcal P^r(x,M)\right\}\right\},$$
observe that $h_u$ is a bounded function and $h_u\uparrow \eta,$ for every $x\in M.$
\begin{align*}
    \mathbb E_{x}\left[\frac{1}{t}\int_0^t f \circ \varphi_s \ \d s \Big\vert \tau > t\right]&=  \frac{1}{t} \frac{
    \displaystyle\int_0^t \mathbb E_{x} \left[f\circ \varphi_s \cdot \mathbbm{1}_{ \tau >t}\right] \d s}{ \mathcal P^t(x,M)}\\
    &= \frac{1}{t} \frac{
    \displaystyle\int_0^t  \mathcal P^s\left(f \cdot   \mathcal P^{t-s}(1)\right)(x) \d s}{ \mathcal P^t(x,M)}\\
    &= \frac{
    \displaystyle \frac{1}{t} \int_0^t e^{s \beta } \mathcal P^s\left(f \cdot e^{(t-s) \beta }  {\mathcal P}^{t-s}(1)\right)(x) \d s}{ e^{\beta  t}\mathcal P^t(x,M)}\\
    &\geq \frac{\eta(x)}{e^{\beta  t}\mathcal P^t(x,M)}
    \displaystyle \frac{1}{\eta(x)}\frac{1}{t} \int_0^{t-u} e^{s\beta  } \mathcal  P^{s}\left(f h_u\right)(x) \d s.
\end{align*}
From Step 1, we get
$$\liminf_{t\to\infty}  \mathbb E_{x}\left[\frac{1}{t}\int_0^t f \circ \varphi_s \ \d s \Big\vert \tau > t\right] \geq \frac{\eta(x) \int_M f h_u \d\mu}{\eta(x)} =  \int_M f h_u \d\mu, \ \text{for every }u\geq 0.$$
Since  $h_u(x)\uparrow \eta(x)$ as $u\to \infty,$ we conclude that  for every $x\in M.$
$$\liminf_{t\to\infty}  \mathbb E_{x}\left[\frac{1}{t}\int_0^t f \circ \varphi_s \ \d s \Big\vert \tau > t\right] \geq \int_M f \eta \d\mu$$

Repeating the same argument to $\|f\|_{\infty} - f,$ we obtain that for every $x\in M,$
$$ \limsup_{t\to\infty}  \mathbb E_{x}\left[\frac{1}{t}\int_0^t f \circ \varphi_s \ \d s \Big\vert \tau > t\right] \leq \int_M f \eta \d\mu.$$
This implies that $\eta(x)\mu(\d x)$ is a quasi-ergodic distribution for $(\theta,\varphi)$ on $M.$ From  $ \mathrm{(H1)} $ we obtain that 
$$\nu(\d x) = \eta(x) \mu(\d x). $$
\begin{step}{3}
We prove $(ii).$
\end{step}
Let $x\in M$ and $A\in\mathcal B(M).$ From Step 1  we obtain that for every $t\geq 0,$
\begin{align*}
    \mu(A) &= \lim_{s\to\infty} \frac{ \mathcal P^{t+s}(x,A)}{\mathcal P^{t+s}(x,M)}\\
    &=\lim_{s\to\infty} 
    \frac{e^{\beta  s}\mathcal P^s \left(e^{\beta  t} \mathcal P^t(\cdot, A)\right)(x)}{e^{\beta  (t+s)}\mathcal P^{t+s}(x,M)}\\
    &= \frac{\eta(x)\int_M  e^{\beta  t} \mathcal P^t(x,A)\mu(\d x)}{\eta(x)} =e^{\beta  t} \int_M\mathcal P^t(x,A) \mu(\d x).
\end{align*}
This proves step 3.
\begin{step}{4}
We prove $(iii).$
\end{step}
From $ \mathrm{(H3)} $ and Lebesgue dominated convergence, 
\begin{align*}
    \mathcal P^t(\eta)(x) &= e^{-\beta  t}\lim_{s\to\infty}\int_M e^{\beta  (s+t)} \mathcal P^s(y,M)\mathcal P^t(x,\d y) \\
    &= e^{-\beta  t}\lim_{s\to\infty} e^{\beta (t+s)}\mathcal P^{t+s}(x,M) = e^{-\beta  t}\eta(x).  
\end{align*}

The proof follows directly from Steps 2--4.
\end{proof}

\end{appendix}

\begin{acks}[Acknowledgments]
We thank Denis Villemonais for an insightful communication about the $Q$-process.
\end{acks}
\begin{funding}
M.C.’s research has been supported by an Imperial College President’s PhD scholarship.  H.C.'s research is funded by a scholarship from the EPSRC DTP in Mathematical Sciences. M.C., H.C. and J.L. have also been supported by the EPSRC Centre for Doctoral Training in Mathematics of Random Systems: Analysis, Modelling and Simulation (EP/S023925/1). J.L. gratefully acknowledges research support from IRCN (Tokyo), GUST (Kuwait), and the London Mathematical Laboratory. D.C. and M.E. thank the DFG SPP 2298 for supporting their research. 
D.C. and M.E. have been additionally supported by Germany's Excellence Strategy -- The Berlin Mathematics Research Center MATH+ (EXC-2046/1, project ID: 390685689).
\end{funding}


\bibliographystyle{imsart-number} 
\bibliography{references.bib}      


\end{document}